\documentclass[12pt]{article}
\usepackage{amsmath,amsthm,amssymb,euscript,verbatim}
\begin{document}
\title
{\bf  Dupin Submanifolds in Lie Sphere Geometry (updated version)}
\author
{Thomas E. Cecil\thanks{Supported by NSF Grant No. DMS 87-06015}
and Shiing-Shen Chern\thanks{Supported by NSF Grant No. DMS 87-01609}}
\maketitle

\begin{abstract}
A hypersurface $M^{n-1}$ in Euclidean space $E^n$ is proper Dupin if the number of distinct principal curvatures
is constant on $M^{n-1}$, and each principal curvature function is constant along each leaf of its principal foliation.  
This paper was originally published in 1989 \cite{CC2}, and it develops a method for the local study of proper Dupin hypersurfaces in the context of Lie sphere geometry using moving frames.
This method has been effective in obtaining several classification theorems of proper Dupin hypersurfaces since that time. This updated version of the paper contains the original exposition together with some remarks by T.Cecil made in 2020 (as indicated in the text)
that describe progress in the field since the time of the original version, as well as some important remaining open problems in the field.
\end{abstract}

\section{Introduction}
\label{sec:1}
Consider a piece of surface immersed in three-dimensional Euclidean space $E^3$.  Its normal lines are the common tangent lines of two surfaces, the focal surfaces.  These focal surfaces may have singularities, and a classical theorem says that if the focal surfaces both degenerate to curves, then the curves are conics, and the surface is a cyclide of Dupin. (See, for example,
\cite[pp. 151--166]{CR7}.)  Equivalently, the cyclides can be characterized as those surfaces in 
$E^3$ whose two distinct principal curvatures are both constant along their corresponding lines of curvature.

The cyclides have been generalized to an interesting class of hypersurfaces in $E^n$, the Dupin hypersurfaces. Initially, a hypersurface $M$ in $E^n$ was said to be Dupin if the number of distinct principal curvatures (or focal points) is constant on $M$ and if each principal curvature is constant along the leaves of its corresponding principal foliation. (See \cite{CR7}, \cite{Th1}, \cite{GH}.) 
More recently, this has been generalized to include cases where the number of distinct principal curvatures is not constant. (See \cite{P4}, \cite{CC1}.)\\

\noindent
{\bf Remark 1.1.} (added by T. Cecil in 2020) In the terminology of Pinkall \cite{P4}, Dupin hypersurfaces for which the number of distinct principal curvatures is constant on $M$ are called {\em proper Dupin}.\\

The study of Dupin hypersurfaces in $E^n$ is naturally situated in the context of Lie sphere geometry, developed by Lie \cite{LS}  as part of his work on contact transformations.  The projectivized cotangent bundle $PT^*E^n$ of $E^n$ has a contact structure.  In fact, if
$x^1,\ldots,x^n$ are the coordinates in $E^n$, the contact structure is defined by the linear differential form 
\begin{displaymath}
dx^n -p_1dx^1-\cdots-p_{n-1} dx^{n-1}.
\end{displaymath}
Lie proved that the pseudo-group of all contact transformations carrying (oriented) hyperspheres in the generalized sense (i.e., including points and oriented hyperplanes) into hyperspheres is a Lie group, called the Lie sphere group, isomorphic to $O(n+1,2)/\pm I$, where $O(n+1,2)$ is the orthogonal group for an indefinite inner product on ${\bf R}^{n+3}$ with signature $(n+1,2)$.  The Lie sphere group contains as a subgroup the M\"{o}bius group of conformal transformations of $E^n$ and, of course, the Euclidean group.  Lie exhibited a bijective correspondence between the set of oriented hyperspheres in $E^n$ and the points on the quadric hypersurface $Q^{n+1}$ in real projective space $P^{n+2}$ given by the equation $\langle x,x \rangle = 0$, where $\langle \ ,\ \rangle$ is the inner product on ${\bf R}^{n+3}$ mentioned above.  The manifold $Q^{n+1}$ contains projective lines but no linear subspaces of $P^{n+2}$ of higher dimension.  The 
1-parameter family of oriented spheres corresponding to the points of a projective line lying on $Q^{n+1}$ consists of all oriented hyperspheres which are in oriented contact at a certain contact element on $E^n$.  Thus, Lie constructed a local diffeomorphism between $PT^*E^n$ and the manifold $\Lambda^{2n-1}$  of projective lines which lie on $Q^{n+1}$.

An immersed submanifold $f:M^k \rightarrow E^n$ naturally induces a Legendre submanifold 
$\lambda: B^{n-1} \rightarrow \Lambda^{2n-1}$, where $B^{n-1}$ is the bundle of unit normal vectors to $f$ (take $B^{n-1} = M^{n-1}$ in the case $k = n-1$).  This Legendre map has similarities with the familiar Gauss map, and like the Gauss map, it can be a powerful tool in the study of submanifolds of Euclidean space.  In particular, the Dupin property for hypersurfaces in $E^n$ is easily formulated in terms of the Legendre map, and it is immediately seen to be invariant under Lie sphere transformations.

The study of Dupin submanifolds has both local and global aspects. Thorbergsson \cite{Th1} showed that a Dupin hypersurface $M$ with $g$ distinct principal curvatures at each point must be taut, i.e., every nondegenerate Euclidean distance function $L_p (x) = |p - x|^2, p \in E^n$, must have the minimum number of critical points on $M$.  Tautness was shown to be invariant Lie sphere transformations in our earlier paper \cite{CC1} (see also \cite[pp. 82--95]{Cec1}).  Using tautness and the work of M\"{u}nzner \cite{Mu}--\cite{Mu2}, Thorbergsson was able to conclude that the number $g$ must be $1,2,3,4$ or 6, as with an isoparametric hypersurface in the sphere $S^n$.
The case $g=1$ is, of course, handled by the well-known classification of umbilic hypersurfaces.  Compact Dupin hypersurfaces with $g=2$ and $g=3$ were classified by Cecil and Ryan \cite{CR2}
and Miyaoka \cite{Mi1}, respectively.  In two recent papers, Miyaoka  \cite{Mi3}, \cite{Mi4}, has made further progress on the classification of compact Dupin hypersurfaces in the cases $g=4$ and $g=6$.  Meanwhile, Grove and Halperin \cite{GH} have determined several important topological invariants of compact Dupin hypersurfaces in the cases $g=4$ and $g=6$. \\

\noindent
{\bf Remark 1.2.} (added by T. Cecil in 2020) For descriptions of more recent progress on the classification of compact proper Dupin hypersurfaces, see \cite[pp. 112--123]{Cec1} and
\cite[pp. 308--322]{CR8}.  For $g \leq 3$, every compact, connected proper Dupin hypersurface with $g$ principal curvatures is equivalent by a Lie sphere transformation to an isoparametric hypersurface in a sphere $S^n.$  This is not true for $g=4$ and $g=6$, however, as shown by examples constructed by Pinkall and Thorbergsson \cite{PT1} for $g=4$,  and by  different examples constructed by Miyaoka and Ozawa \cite{MO} for $g=4$ and $g=6$.  These examples do not have constant Lie curvatures (cross-ratios of the principal curvatures, see Miyaoka \cite{Mi3}), so they can't be Lie equivalent to an isoparametric hypersurface, which must have constant Lie curvatures.
The classification of compact proper Dupin hypersurfaces in the cases $g=4$ and  $g=6$ is still an open problem.\\

In this paper, we study Dupin hypersurfaces in the setting of Lie sphere geometry using local techniques.  In Section 2, we give a brief introduction to Lie sphere geometry.  In Section 3, we introduce the basic differential geometric notions: the Legendre map and the Dupin property.  The case of $E^3$ is handled in Section 4, where we handle the case of $g=2$ distinct focal points for $E^n$.  This was first done for $n>3$ by Pinkall \cite{P4}.  Our main contribution lies in Section 5, where we treat the case $E^4$ by the method of moving frames.  This case was also studied by Pinkall \cite{P3}, but our treatment seems to be more direct and differs from his in several points.  It is our hope that this method will provide a framework and give some direction for the study of Dupin hypersurfaces in $E^n$ for $n>4$.

\section{Lie Sphere Geometry}
\label{sec:2}

We first present a brief outline of the main ideas in Lie's geometry of spheres in ${\bf R}^n$.  This is given in more detail in Lie's original treatment \cite{LS}, in the book of Blaschke \cite{Bl}, and in our paper \cite{CC1} (see also the book \cite{Cec1}).

The basic construction in Lie sphere geometry associates each oriented sphere, oriented plane and point sphere in ${\bf R}^n \cup \{ \infty \} = S^n$ with a point on the quadric $Q^{n+1}$ in projective space $P^{n+2}$ given in homogeneous coordinates $(x_1,\ldots,x_{n+3})$ by the equation
\begin{equation}
\label{eq:2.1}
\langle x,x \rangle = - x_1^2 + x_2^2 + \cdots + x_{n+2}^2 - x_{n+3}^2 = 0.
\end{equation}
We will denote real $(n+3)$-space endowed with the metric (\ref{eq:2.1}) by ${\bf R}^{n+3}_2$.

We can designate the orientation of a sphere in ${\bf R}^n$ by assigning a plus or a minus sign to its radius.  Positive radius corresponds to the orientation determined by the field of inward normals to the sphere, while a negative radius corresponds to the orientation determined by the outward normal.  (See Remark 2.1 below.) A plane in ${\bf R}^n$ is a sphere which goes through the point 
$\infty$.  The orientation of the plane can be associated with a choice of unit normal $N$.  The specific correspondence between the points of $Q^{n+1}$ and the set of oriented spheres, oriented planes and points in ${\bf R}^n \cup \{ \infty \}$ is then given as follows:

\begin{equation}
\label{eq:2.2}
\begin{array}{cc}
{\rm {\bf Euclidean}} & {\rm {\bf Lie}} \\
 & \\
\mbox{{\rm points:} }u \in {\bf R}^n & \left[ \left( \frac{1+u \cdot u}{2},\frac{1-u \cdot u}{2},u,0 \right) 
\right]  \\
 & \\
\infty & [(1,-1,0,0)]\\
 & \\
\mbox{{\rm spheres: center} $p$, {\rm signed radius} $r$} & \left[ \left( \frac{1+p \cdot p - r^2}{2}, 
\frac{1- p \cdot p +r^2} {2},p, r \right) \right] \\
 & \\ 
\mbox{{\rm planes:} $u \cdot N = h$, {\rm unit normal} $N$} & [(h,-h,N,1)] 
\end{array}
\end{equation}

\vspace{.25in}

\noindent
Here the square brackets denote the point in projective space $P^{n+2}$ given by the homogeneous coordinates in the round brackets, and $u \cdot u$ is the standard Euclidean dot product in ${\bf R}^n$.  

From (\ref{eq:2.2}), we see that the point spheres correspond to the points in the intersection of $Q^{n+1}$ with the hyperplane in $P^{n+2}$ given by the equation $x_{n+3} = 0.$  The manifold of point spheres is called {\em M\"{o}bius space}.

A fundamental notion in Lie sphere geometry is that of oriented contact of spheres.  Two oriented spheres $S_1$ and $S_2$ are in {\em oriented contact} if they are tangent and their orientations agree at the point of tangency.  If $p_1$ and $p_2$ are the respective centers of $S_1$ and $S_2$, and $r_1$ and $r_2$ are the respective signed radii, then the condition of oriented contact can be expressed analytically by
\begin{equation}
\label{eq:2.3}
|p_1 - p_2| = |r_1 - r_2|.
\end{equation}
If $S_1$ and $S_2$ are represented by $[k_1]$ and $[k_2]$ as in (\ref{eq:2.2}), then (\ref{eq:2.3}) is equivalent to the condition
\begin{equation}
\label{eq:2.4}
\langle k_1 , k_2 \rangle = 0.
\end{equation}

In the case where $S_1$ and/or $S_2$ is a plane or a point in ${\bf R}^n$, oriented contact has the logical meaning.  That is, a sphere $S$ and a plane $\pi$ are in oriented contact if $\pi$ is tangent to $S$ and their orientations agree at the point of contact.  Two oriented planes are in oriented contact if their unit normals are the same.  They are in oriented contact at the point $\infty$.  A point sphere is in oriented contact with a sphere or plane $S$ if it lies on $S$, and $\infty$ is in oriented contact with each plane.  In each case, the analytic condition for oriented contact is equivalent to (\ref{eq:2.4}) when the two ``spheres'' in question are represented as in (\ref{eq:2.2}).\\

\noindent
{\bf Remark 2.1.} (added by T. Cecil in 2020) In the case of a sphere $[k_1]$ and a plane $[k_2]$ as in (\ref{eq:2.2}), the equation (\ref{eq:2.4}) is equivalent to $p \cdot N = h + r$. In order to make this correspond to the geometric definition of oriented contact stated in the paragraph above, one must adopt the convention that the inward normal orientation of the sphere corresponds to a positive signed radius. \\

Because of the signature of the metric (\ref{eq:2.1}), the quadric $Q^{n+1}$ contains lines in $P^{n+2}$ but no linear subspaces of higher dimension.  A line in $Q^{n+1}$  is determined by two points $[x], [y]$ in $Q^{n+1}$ satisfying $\langle x,y \rangle = 0$.  The lines on $Q^{n+1}$ form a manifold of dimension $2n-1$, to be denoted by $\Lambda^{2n-1}$.  In ${\bf R}^n$, a line on $Q^{n+1}$ corresponds to a 1-parameter family of oriented spheres such that any two of the oriented spheres are in oriented contact, i.e., all the oriented spheres tangent to an oriented plane at a given point, i.e., an oriented contact element.  Of course, a contact element can also be represented by an element of $T_1S^n$, the bundle of unit tangent vectors to the Euclidean sphere $S^n$ in $E^{n+1}$ with its usual metric.  This is the starting for Pinkall's \cite{P4} considerations of Lie geometry.

A {\em Lie sphere transformation} is a projective transformation of $P^{n+2}$ which takes $Q^{n+1}$ to itself.  Since a projective transformation takes lines to lines, a Lie sphere transformation preserves oriented contact of spheres.  The group $G$ of Lie sphere transformations is isomorphic to $O(n+1,2)/\{\pm I\}$, where $O(n+1,2)$ is the group of orthogonal transformations for the inner product (\ref{eq:2.1}).  The group of M\"{o}bius transformations is isomorphic to $O(n+1,1)/\{\pm I\}$.

\section{Legendre Submanifolds}
\label{sec:3}

Here we recall the concept of a Legendre submanifold of the contact manifold $\Lambda^{2n-1} (= \Lambda)$  using the notation of \cite{CC1}.  In this section, the ranges of the indices are as follows:
\begin{equation}
\label{eq:3.1}
1 \leq A, B, C \leq n+3, \quad 3 \leq i,j,k \leq n+1.
\end{equation}
Instead of using an orthornormal frame for the metric $\langle \ ,\  \rangle$ defined by (\ref{eq:2.1}), it is useful to consider a {\em Lie frame}, that is an ordered set of vectors $Y_A$ in 
${\bf R}^{n+3}_2$ satisfying
\begin{equation}
\label{eq:3.2}
\langle Y_A, Y_B \rangle = g_{AB},
\end{equation}
with
\begin{equation}
\label{eq:3.3}
[g_{AB}] = \left[ \begin{array}{ccc}J&0&0\\0&I_{n-1}&0\\ 0&0&J\end{array}\right],
\end{equation}
where $I_{n-1}$ is the $(n-1) \times (n-1)$ identity matrix and
\begin{equation}
\label{eq:3.4}
J = \left[ \begin{array}{cc}0&1\\1&0\end{array}\right].
\end{equation}
The space of all Lie frames can be identified with the orthogonal group $O(n+1,2)$, of which the Lie sphere group, being isomorphic to $O(n+1,2)/\{\pm I\}$, is a quotient group.  In this space, we introduce the Maurer-Cartan forms
\begin{equation}
\label{eq:3.5}
dY_A = \sum \omega_A^B Y_B.
\end{equation}
Through differentiation of (\ref{eq:3.2}), we show that the following matrix of 1-forms is skew-symmetric
\vspace{.2in}
\begin{equation}
\label{eq:3.6}
[\omega_{AB}] = \left[ \begin{array}{ccccc}\omega_1^2&\omega_1^1&\omega_1^i&\omega_1^{n+3}&\omega_1^{n+2}\\
\omega_2^2&\omega_2^1&\omega_2^i&\omega_2^{n+3}&\omega_2^{n+2}\\
\omega_j^2&\omega_j^1&\omega_j^i&\omega_j^{n+3}&\omega_j^{n+2}\\
\omega_{n+2}^2&\omega_{n+2}^1&\omega_{n+2}^i&\omega_{n+2}^{n+3}&\omega_{n+2}^{n+2}\\
\omega_{n+3}^2&\omega_{n+3}^1&\omega_{n+3}^i&\omega_{n+3}^{n+3}&\omega_{n+3}^{n+2}\\
\end{array}\right].
\end{equation}
\vspace{.2in}
Next by taking the exterior derivative of (\ref{eq:3.5}), we get the Maurer-Cartan equations
\begin{equation}
\label{eq:3.7}
d\omega_A^B = \sum_C \omega_A^C \wedge \omega_C^B.  
\end{equation}
In \cite{CC1}, we then show that the form
\begin{displaymath}
\omega_1^{n+2} = \langle dY_1, Y_{n+3} \rangle
\end{displaymath}
gives a contact structure on the manifold $\Lambda$.

Let $B^{n-1} (=B)$ be an $(n-1)$-dimensional smooth manifold.  A {\em Legendre map} is a smooth map $\lambda:B \rightarrow \Lambda$ which annihilates the contact form on $\Lambda$, i.e., 
$\lambda^* \omega_1^{n+2}  = 0$ on $B$.  All of our calculations are local in nature.  We use the method of moving frames  and consider a smooth family of Lie frames $Y_A$ on an open subset $U$ of $B$, with the line $\lambda(b)$ given by $[Y_1 (b), Y_{n+3} (b)]$ for each $b \in U$.  The Legendre map $\lambda$ is called a {\em Legendre submanifold}  if for a generic choice of $Y_1$ the forms 
$\omega_1^i,\  3 \leq i \leq n+1$, are linearly independent, i.e.,
\begin{equation}
\label{eq:3.8}
\wedge \ \omega_1^i \neq 0 \  {\rm on}\ U.  
\end{equation}
Here and later we pull back the structure forms to $B^{n-1}$ and omit the symbols of such pull-backs for simplicity.  Note that the Legendre condition is just
\begin{equation}
\label{eq:3.9}
\omega_1^{n+2} = 0.  
\end{equation}
We now assume that our choice of $Y_1$ satisfies (\ref{eq:3.8}).  By exterior differentiation of (\ref{eq:3.9}) and using (\ref{eq:3.6}), we get 
\begin{equation}
\label{eq:3.10}
\sum \omega_1^i \wedge \omega_{n+3}^i = 0. 
\end{equation}
Hence by Cartan's Lemma and  (\ref{eq:3.8}), we have
\begin{equation}
\label{eq:3.11}
\omega_{n+3}^i = \sum h_{ij} \omega_1^j, \ {\rm with}\ h_{ij} = h_{ji}. 
\end{equation}
The quadratic differential form
\begin{equation}
\label{eq:3.12}
II(Y_1) = \sum_{i,j} h_{ij} \omega_1^i \omega_1^j,
\end{equation}
defined up to a non-zero factor and depending on the choice of $Y_1$, is called the {\em second fundamental form}.

This form can be related to the well-known Euclidean second fundamental form in the following way.  Let $e_{n+3}$ be any unit timelike vector in ${\bf R}^{n+3}_2$.  For each $b \in U$, let $Y_1(b)$ be the point of intersection of the line $\lambda (b)$ with the hyperplane $e_{n+3}^\perp$.
$Y_1$ represents the locus of point spheres in the M\"{o}bius space $Q^{n+1} \cap e_{n+3}^\perp$, and we call $Y_1$ the {\em M\"{o}bius  projection of} $\lambda$ determined by $e_{n+3}$.
Let $e_1$ and $e_2$ be unit timelike, respectively spacelike, vectors orthogonal to $e_{n+3}$ and to each other, chosen so that $Y_1$ is not the point at infinity $[e_1 - e_2]$ for any $b \in U$.  We can represent $Y_1$ by the vector
\begin{equation}
\label{eq:3.13}
Y_1 = \frac{1 + f \cdot f}{2} e_1 + \frac{1 - f \cdot f}{2} e_2 + f, 
\end{equation}
as in (\ref{eq:2.2}), where $f(b)$ lies in the space ${\bf R}^n$ of vectors orthogonal to $e_1, e_2$ and $e_{n+3}$.  We will call the map $f:B \rightarrow {\bf R}^n$ the {\em Euclidean projection} of $\lambda$ determined by the ordered triple $e_1, e_2, e_{n+3}$.  The regularity condition 
(\ref{eq:3.8}) is equivalent to the condition that $f$ be an immersion on $U$ into ${\bf R}^n$.  For each $b \in U$, let $Y_{n+3} (b)$ be the intersection of $\lambda (b)$ with the orthogonal complement of the lightlike vector $e_1 - e_2$. $Y_{n+3}$ is distinct from $Y_1$ and thus
$\langle Y_{n+3} , e_{n+3} \rangle \neq 0$.  So we can represent $Y_{n+3}$ by a vector of the form
\begin{equation}
\label{eq:3.14}
Y_{n+3} = h (e_1 - e_2) + \xi +e_{n+3},
\end{equation}
where $\xi:U \rightarrow {\bf R}^n$ has unit length and $h$ is a smooth function on $U$.  Thus, according to  (\ref{eq:2.2}), $Y_{n+3} (b)$ represents the plane in the pencil of oriented spheres in 
${\bf R}^n$ corresponding to the line $\lambda (b)$ on $Q^{n+1}$.  Note that the condition 
$\langle Y_1, Y_{n+3}  \rangle = 0$ is equivalent to $h = f \cdot \xi$, while the Legendre condition 
$\langle dY_1, Y_{n+3}  \rangle = 0$ is the same as the Euclidean condition 
\begin{equation}
\label{eq:3.15}
\xi \cdot df = 0.
\end{equation}
Thus, $\xi$ is a field of unit normals to the immersion $f$ on $U$.  Since $f$ is an immersion, we can choose the Lie frame vectors $Y_3,\ldots,Y_{n+1}$ to satisfy
\begin{equation}
\label{eq:3.16}
Y_i = dY_1 (X_i) = (f \cdot df(X_i))(e_1-e_2) + df(X_i),\ 3 \leq i \leq n+1,
\end{equation}
for tangent vector fields $X_3,\ldots,X_{n+1}$ on $U$.  Then, we have 
\begin{equation}
\label{eq:3.17}
\omega_1^i (X_j) = \langle dY_1 (X_j), Y_i \rangle = \langle Y_j, Y_i \rangle = \delta_{ij}.
\end{equation}
Now using (\ref{eq:3.14}) and (\ref{eq:3.16}), we compute
\begin{equation}
\label{eq:3.18}
\omega_{n+3}^i (X_j) = \langle dY_{n+3} (X_j), Y_i \rangle =  d\xi (X_j) \cdot df(X_i)  = 
- df(AX_j) \cdot df(X_i) = -A_{ij},
\end{equation}
where $A = [A_{ij}]$ is the Euclidean shape operator (second fundamental form) of the immersion $f$.  But by (\ref{eq:3.11}) and (\ref{eq:3.17}), we have
\begin{displaymath}
\omega_{n+3}^i (X_j) = \sum h_{ik} \omega_1^k (X_j) = h_{ij}.
\end{displaymath}
Hence $h_{ij} = -A_{ij}$, and $[h_{ij}]$ is just the negative of the Euclidean shape operator $A$ of $f$.\\

\noindent
{\bf Remark 3.1}  The discussion above demonstrates how an immersion 
$f:B^{n-1} \rightarrow {\bf R}^n$ with field of unit normals $\xi$ induces a Legendre submanifold 
$\lambda:B^{n-1} \rightarrow \Lambda$ defined by $\lambda (b) = [Y_1(b), Y_{n+3} (b)]$, for 
$Y_1, Y_{n+3}$ as in (\ref{eq:3.13}) and
(\ref{eq:3.14}), called the {\em Legendre lift}
of the immersion $f$ with field of unit normals $\xi$.
Further, an immersed submanifold $f:M^k \rightarrow {\bf R}^n$ of codimension greater than one gives rise to a Legendre submanifold $\lambda:B^{n-1} \rightarrow \Lambda$, where $B^{n-1}$ is the bundle of unit normals to $f$ in ${\bf R}^n$.  As in the case of codimension one, $\lambda (b)$ is defined to be the line on $Q^{n+1}$ corresponding to the oriented contact element determined by the unit vector $b$ normal to $f$ at the point $x = \pi (b)$, where $\pi$ is the bundle projection from $B^{n-1}$ to $M^k$.  (In that case, $\lambda$ is called the Legendre lift of the submanifold $f$.)\\

\noindent
{\bf Remark 3.2.} (added by T. Cecil in 2020) 
In a similar way, an immersion $\varphi:B^{n-1} \rightarrow S^n$ with field of unit normals 
$\eta:B^{n-1} \rightarrow S^n$ in $S^n$
induces a Legendre submanifold $\mu:B^{n-1} \rightarrow \Lambda$ defined by
$\mu (b) = [Z_1(b), Z_{n+3} (b)]$, where
\begin{equation}
\label{eq:18a}
Z_1(b) = e_1 + \varphi (b), \quad Z_{n+3}(b) = \eta (b) + e_{n+3},
\end{equation} 
and $S^n$ is identified with the unit sphere in the space
${\bf R}^{n+1}$ of vectors orthogonal to $e_1$ and $e_{n+3}$.  
The map $\mu$ is called the {\em Legendre lift} of the immersion $\varphi$ with field of unit normals 
$\eta$.  As in the Euclidean case, an immersed submanifold $\varphi:M^k \rightarrow S^n$ of codimension greater than one gives rise to a Legendre submanifold 
$\mu:B^{n-1} \rightarrow \Lambda$, where $B^{n-1}$ is the bundle of unit normals to $\varphi$ in 
$S^n$.  As in the case of codimension one, $\mu (b)$ is defined to be the line on $Q^{n+1}$ corresponding to the oriented contact element determined by the unit vector $b$ normal to 
$\varphi$ at the point $x = \pi (b)$, where $\pi$ is the bundle projection from $B^{n-1}$ to $M^k$.
In that case, $\mu$ is called the Legendre lift of the submanifold $\varphi$.\\

As one would expect, the eigenvalues of the second fundamental form have geometric significance.  Consider a curve $\gamma (t)$ on $B$.  The set of points in $Q^{n+1}$ lying on the lines 
$\lambda (\gamma (t))$ forms a ruled surface in $Q^{n+1}$.  We look for the conditions that this ruled surface be developable, i.e., consist of tangent lines to a curve in $Q^{n+1}$.  Let 
$rY_1 + Y_{n+3}$ be the point of contact.  We have by (\ref{eq:3.5}) and (\ref{eq:3.6})
\begin{equation}
\label{eq:3.19}
d(rY_1 + Y_{n+3} ) \equiv \sum_i (r \omega_1^i + \omega_{n+3}^i) Y_i, \bmod \ Y_1, Y_{n+3}.
\end{equation}
Thus, the lines $\lambda (\gamma (t))$ form a developable if and only if the tangent direction of 
$\gamma (t)$ is a common solution to the equations
\begin{equation}
\label{eq:3.20}
\sum_j (r \delta_{ij} + h_{ij} )\   \omega_1^j = 0, \ 3\leq i \leq n+1.
\end{equation}
In particular, $r$ must be a root of the equation
\begin{equation}
\label{eq:3.21}
\det (r \delta_{ij} + h_{ij} ) = 0.
\end{equation}
By (\ref{eq:3.11}) the roots of (\ref{eq:3.21}) are all real.  Denote them by $r_3,\ldots,r_{n+1}$.  The points $r_i Y_1 + Y_{n+3}, \ 3 \leq i \leq n+1$, are called the {\em focal points} or {\em curvature spheres} (Pinkall \cite{P4}) on $\lambda (b)$.  If $Y_1$ and $Y_{n+3}$ correspond to an immersion $f:U \rightarrow  {\bf R}^n$ as in (\ref{eq:3.13}) and (\ref{eq:3.14}), then these focal points on 
$\lambda (b)$ correspond by (\ref{eq:2.2}) to oriented spheres in ${\bf R}^n$ tangent to $f$ at $f(b)$ and centered at the Euclidean focal points of $f$.  These spheres are called the curvature spheres of $f$ and the $r_i$ are just the principal curvatures of $f$, i.e., eigenvalues of the shape operator $A$.  

If $r$ is a root of multiplicity $m$, then the equations (\ref{eq:3.20}) define an $m$-dimensional subspace $T_r$ of $T_b B$, the tangent space to $B$ at the point $b$.  The space $T_r$ is called a {\em principal space} of $T_b B$, the latter being decomposed into a direct sum of its principal spaces.  Vectors in $T_r$ are  called {\em principal vectors} corresponding to the focal point 
$rY_1 + Y_{n+3}$.  Of course, if $Y_1$ and $Y_{n+3}$ correspond to an immersion 
$f:U \rightarrow {\bf R}^n$ as in (\ref{eq:3.13}) and (\ref{eq:3.14}), then these principal vectors are the same as the Euclidean principal vectors for $f$ corresponding to the principal curvature $r$.

With a change of frame of the form
\begin{equation}
\label{eq:3.22}
Y_i^* =  \sum_i c_i^j Y_j, \  3\leq i \leq n+1,
\end{equation}
where $[c_i^j]$ is an $(n-1) \times (n-1)$ orthogonal matrix, we can diagonalize $[h_{ij}]$ so that in the new frame, equation (\ref{eq:3.11}) has the form
\begin{equation}
\label{eq:3.23}
\omega_{n+3}^i = - r_i \omega_1^i, \ 3\leq i \leq n+1.
\end{equation}
Note that none of the functions $r_i$ is ever infinity on $U$ because of the assumption that
(\ref{eq:3.8}) holds, i.e., $Y_1$ is not a focal point.  By associating $Y_1$ to a Euclidean immersion $f$ as in (\ref{eq:3.13}), we can apply results from Euclidean geometry to our situation.  In particular, it follows from a result of Singley \cite{Sin} on Euclidean shape operators that there is a dense open subset of $B$ on which the number $g(b)$ of distinct focal points on $\lambda (b)$ is locally constant.  We will work exclusively on open subsets $U$ of $B$ on which $g$ is constant.  In that case, each eigenvalue function $r:U \rightarrow {\bf R}$ is smooth (see Nomizu \cite{Nom2}), and its corresponding principal distribution is a smooth $m$-dimensional foliation, where $m$ is the multiplicity of $r$ (see \cite[p. 139]{CR7}).  Thus, on $U$ we can find smooth vector fields $X_3,\ldots,X_{n+1}$ dual to smooth 1-forms $\omega_1^3,\ldots,\omega_1^{n+1}$, respectively, such that each $X_i$ is principal for the smooth focal point map $r_i Y_1 + Y_{n+3}$ on $U$.  If
$r Y_1 + Y_{n+3}$ is a smooth focal point map of multiplicity $m$ on $U$, then we can assume that
\begin{equation}
\label{eq:3.24}
r_3 = \cdots = r_{m+2} = r.
\end{equation}
By a different choice of the point at infinity, i.e., $e_1$ and $e_2$, if necessary, we can also assume that the function $r$ is never zero on $U$, i.e., $Y_{n+3}$ is not a focal point on $U$.

We now want to consider a Lie frame $Y_A^*$ for which $Y_1^*$ is a smooth focal point map of multiplicity $m$ on $U$.  Specifically, we make the change of frame
\begin{eqnarray}
\label{eq:3.25}
Y_1^*& =& r Y_1 + Y_{n+3} \nonumber \\
Y_2^* &=& (1/r) Y_2 \nonumber \\
Y_{n+2}^* &=& Y_{n+2} - (1/r) Y_2 \\
Y_{n+3}^*&=&Y_{n+3} \nonumber \\
Y_i^*&=& Y_i, \  3\leq i \leq n+1. \nonumber
\end{eqnarray}
We denote the Maurer-Cartan forms in this frame by $\theta_A^B$.  Note that
\begin{equation}
\label{eq:3.26}
dY_1^* = d(rY_1 + Y_{n+3})  = (dr) Y_1 + r dY_1 +dY_{n+3} = \sum \theta_1^A Y_A^*.
\end{equation}
By examining the coefficient of $Y_i^* = Y_i$ in (\ref{eq:3.26}) , we see from (\ref{eq:3.23}) that
\begin{equation}
\label{eq:3.27}
\theta_1^i = r \omega_1^i + \omega_{n+3}^i = (r - r_i)\  \omega_1^i, \ 3\leq i \leq n+1.
\end{equation}
From (\ref{eq:3.24}) and (\ref{eq:3.27}), we see that 
\begin{equation}
\label{eq:3.28}
\theta_1^a= 0,  \ 3\leq a \leq m+2.
\end{equation}
This equation characterizes the condition that a focal point map $Y_1^*$ have constant multiplicity $m$ on $U$.

We now introduce the concept of a Dupin submanifold and then see what further restrictions it allows us to place on the structure forms.

A connected submanifold $N \subset B$ is called a {\em curvature submanifold} if its tangent space is everywhere a principal space.  The Legendre submanifold is called {\em Dupin} if for every curvature submanifold $N \subset B$, the lines $\lambda (b), b \in N$, all pass through a fixed point, i.e., each focal point map is constant along its curvature submanifolds.
This definition of ``Dupin'' is the same as that of Pinkall \cite{P4}.  It is weaker than the definition of Dupin for Euclidean hypersurfaces used by Thorbergsson \cite{Th1}, Miyaoka \cite{Mi1}, Grove-Halperin \cite{GH} and Cecil-Ryan \cite{CR7}, all of whom assumed that the number $g$ of distinct curvature spheres is constant on $B$.  As we noted above, $g$ is locally constant on a dense open subset of any Legendre submanifold, but $g$ is not necessarily constant on the whole of a Dupin submanifold, as the example of the Legendre submanifold induced from a tube $B^3$ over a torus $T^2 \subset {\bf R}^3 \subset {\bf R}^4$ demonstrates (see Pinkall \cite{P3}, Cecil-Ryan 
\cite[p. 188]{CR7}).

It is easy to see that the Dupin property is invariant under Lie sphere transformations as follows.  Suppose that $\lambda:B \rightarrow \Lambda^{2n-1}$ is a Legendre submanifold and that $\alpha$ is a Lie transformation.  The map $\alpha \lambda:B \rightarrow \Lambda^{2n-1}$ is also a Legendre submanifold with $\alpha \lambda (b) = [\alpha Y_1 (b), \alpha Y_{n+3} (b)]$.  Furthermore, if $k = r Y_1 + Y_{n+3}$ is a curvature sphere of $\lambda$ at a point $b \in B$, then since $\alpha$ is a linear transformation, $\alpha k$ is a curvature sphere of $\alpha \lambda$ at $b$ with the same principal space as $k$.  Thus $\lambda$ and $\alpha \lambda$ have the same curvature submanifolds on $B$, and the Dupin property clearly holds for $\lambda$ if and only if it holds for $\alpha \lambda$.

We now return to the calculation that led to equation (\ref{eq:3.28}).  We have that $Y_1^* = r Y_1 + Y_{n+3}$ is a smooth focal point map of multiplicity $m$ on the open set $U$ and its corresponding principal space is spanned by the vector fields $X_3,\ldots,X_{m+2}$.  The Dupin condition that $Y_1^*$ be constant along the leaves of $T_r$ is simply
\begin{equation}
\label{eq:3.29}
dY_1^* (X_a) \equiv 0,\  \bmod \ Y_1^*, \ 3\leq a \leq m+2.
\end{equation}
From (\ref{eq:3.17}), (\ref{eq:3.27}) and (\ref{eq:3.28}), we have that
\begin{equation}
\label{eq:3.30}
dY_1^* (X_a) = \theta_1^1 (X_a) Y_1 + \theta_1^{n+3} (X_a) Y_{n+3}, \ 3\leq a \leq m+2.
\end{equation}
Comparing (\ref{eq:3.29}) and (\ref{eq:3.30}), we see that
\begin{equation}
\label{eq:3.31}
\theta_1^{n+3} (X_a) = 0, \ 3\leq a \leq m+2.
\end{equation}
We now show that we can make one more change of frame and make $\theta_1^{n+3} = 0$.  We can write the form $\theta_1^{n+3}$ in terms of the basis $\omega_1^3,\ldots, \omega_1^{n+1}$ as
\begin{displaymath}
\theta_1^{n+3} = \sum a_i \omega_1^i.
\end{displaymath}
From (\ref{eq:3.31}), we see that we actually have
\begin{equation}
\label{eq:3.32}
\theta_1^{n+3} = \sum_{b=m+3}^{n+1} a_b \omega_1^b .
\end{equation}
Using (\ref{eq:3.17}), (\ref{eq:3.27}) and (\ref{eq:3.32}), we compute for all $m+3 \leq b \leq n+1$,
\begin{eqnarray}
\label{eq:3.33}
dY_1^* (X_b) & = & \theta_1^1(X_b) Y_1^* + \theta_1^b(X_b) Y_b + \theta_1^{n+3}(X_b) Y_{n+3}  \nonumber \\
& = & \theta_1^1(X_b) Y_1^* +(r - r_b) Y_b + a_b Y_{n+3} \\
& = & \theta_1^1(X_b) Y_1^* +(r - r_b)(Y_b +(a_b/(r - r_b)) Y_{n+3}).\nonumber
\end{eqnarray}
We now make the change of Lie frame,
\begin{eqnarray}
\label{eq:3.34}
\overline{Y}_1 & = & Y_1^*, \quad \overline{Y}_2 = Y_2^*, \nonumber \\
 \overline{Y}_a &=& Y_a, \quad 3 \leq a \leq m+2, \nonumber \\
\overline{Y}_b & = & Y_b + (a_b/(r - r_b))Y_{n+3}, \quad m+3 \leq b \leq n+1, \\
\overline{Y}_{n+2} & = & - \sum_{b=m+3}^{n+1}(a_b/(r - r_b)) Y_b + Y_{n+2} - \frac{1}{2}
\sum_{b=m+3}^{n+1}(a_b/(r - r_b))^2 Y_{n+3}, \nonumber\\
\overline{Y}_{n+3}  &=& Y_{n+3}. \nonumber
\end{eqnarray}
Let $\alpha_A^B$ be the Maurer--Cartan forms for this new frame.  We still have
\begin{equation}
\label{eq:3.35}
\alpha_1^a = \langle d\overline{Y}_1, \overline{Y}_a \rangle = \langle dY_1^*, Y_a \rangle = \theta_1^a = 0, \quad 3 \leq a \leq m+2.
\end{equation}
Furthermore, since $\overline{Y}_1 = Y_1^*$, the Dupin condition (\ref{eq:3.29}) still yields
\begin{equation}
\label{eq:3.36}
\alpha_1^{n+3} (X_a) = 0, \quad 3 \leq a \leq m+2.
\end{equation}
Finally, for $m+3 \leq b \leq n+1$,  we have from (\ref{eq:3.33}) and (\ref{eq:3.34})
\begin{equation}
\label{eq:3.37}
\alpha_1^{n+3} (X_b) = \langle d\overline{Y}_1(X_b), \overline{Y}_{n+2} \rangle = \langle \theta_1^1(X_b) \overline{Y}_1 + (r - r_b)\overline{Y}_b, \overline{Y}_{n+2} \rangle 
= 0.
\end{equation}
From (\ref{eq:3.36}) and (\ref{eq:3.37}), we conclude that
\begin{equation}
\label{eq:3.38}
\alpha_1^{n+3} = 0.
\end{equation}
Thus, our main result of this section is that the assumption that the focal point map
$\overline{Y}_1 = r Y_1 + Y_{n+3}$ has constant multiplicity $m$ and is constant along the leaves of its principal foliation $T_r$ allows us to produce a Lie frame $\overline{Y}_A$ whose structure forms satisfy
\begin{eqnarray}
\label{eq:3.39}
\alpha_1^a &=& 0, \ 3 \leq a \leq m+2,\\
\alpha_1^{n+3} &=& 0. \nonumber
\end{eqnarray}

\section{Cyclides of Dupin}
\label{sec:4}

Dupin \cite[p. 200]{D} initiated the study of this subject in 1822 when he defined a cyclide to be a surface $M^2$ in $E^3$ which is the envelope of the family of spheres tangent to three fixed spheres.  This was shown to be equivalent requiring that both sheets of the focal set of $M^2$ in $E^3$ degenerate into curves.  Then $M^2$ is the envelope of each of the two families of curvature spheres.  The key step in the classical Euclidean proof (see, for example, Eisenhart 
\cite[pp. 312--314]{Eisenhart} or Cecil-Ryan 
\cite[pp. 151--166]{CR7}) is to show that the two focal curves are a pair of so-called ``focal conics'' in $E^3$, i.e., an ellipse and hyperbola in mutually orthogonal planes such that the vertices of the ellipse are the foci of the hyperbola and vice-versa, or a pair of parabolas in orthogonal planes such that the vertex of each is the focus of the other.  This classical proof is local, i.e., one needs only a small piece of the surface to determine the focal conics and reconstruct the whole cyclide.  Of course, envelopes of families of spheres can have singularities and some of the cyclides have one or two singular points in $E^3$.  It turns out, however, that all of the different forms of cyclides in Euclidean space induce Legendre submanifolds which are locally Lie equivalent.  In other words, they are all various Euclidean projections of one Legendre submanifold.  Pinkall \cite{P4} generalized this result to higher dimensional Dupin submanifolds.  He defined a {\em cyclide of characteristic} $(p,q)$ to be a Dupin submanifold with the property that at each point it has exactly two distinct focal points with respective multiplicities $p$ and $q$.  He then proved the following.\\

\noindent
{\bf Theorem 4.1.} (Pinkall \cite{P4}): {\it (a) Every connected cyclide of Dupin is contained in a unique compact and connected cyclide of Dupin.}\\
{\it (b) Any two cyclides of the same characteristic are locally Lie equivalent, each being Lie equivalent to an open subset of a standard product of spheres in $S^n$.}\\

In this section, we give a proof of Pinkall's result using the method of moving frames.  Let 
$\lambda:B \rightarrow \Lambda$ be the Dupin cyclide.  The main step in the proof of Theorem 4.1 is to show that the two focal point maps $k_1$ and $k_2$ from $B$ to $Q^{n+1}$ are such that the image $k_1 (B)$ lies in the intersection of $Q^{n+1}$ with a $(p+1)$-dimensional subspace $E$ of $P^{n+2}$ while $k_2 (B)$ lies in the intersection of the $(q+1)$-dimensional subspace $E^\perp$ with $Q^{n+1}$.  This generalizes the key step in the classical Euclidean proof that the two focal curves are focal conics.  Once this fact has been established for $k_1$ and $k_2$, it is relatively easy to complete the proof of the Theorem. 

We begin the proof by taking advantage of the results of the previous section.  As we showed in (\ref{eq:3.39}), on any neighborhood $U$ in $B$, we can find a local Lie frame, which we now denote by $Y_A$, whose Maurer-Cartan forms, now denoted $\omega_A^B$, satisfy
\begin{eqnarray}
\label{eq:4.1}
\omega_1^a &=& 0, \ 3 \leq a \leq p+2,\\
\omega_1^{n+3} &=& 0. \nonumber
\end{eqnarray}
In this frame, $Y_1$ is a focal point map of multiplicity $p$ from $U$ to $Q^{n+1}$.  By the hypotheses of Theorem 4.1, there is one other focal point of multiplicity $q = n-1-p$ at each point of $B$.  By repeating the procedure used in constructing the frame $Y_A$, we can construct a new frame $\overline{Y}_A$ which has as $\overline{Y}_{n+3}$ the other focal point map 
$sY_1 + Y_{n+3}$, where $s$ is a root of (\ref{eq:3.21}) of multiplicity $q$.  The principal space $T_s$ is the span of the vectors $X_{p+3},\ldots,X_{n+1}$ in the notation of the previous section.
The fact that $\overline{Y}_{n+3}$ is a focal point yields
\begin{equation}
\label{eq:4.2}
\overline{\omega}_{n+3}^b= 0,  \ p+3\leq b \leq n+1,
\end{equation}
in analogy to (\ref{eq:3.28}).  The Dupin condition analogous to (\ref{eq:3.29}) is
\begin{equation}
\label{eq:4.3}
d\overline{Y}_{n+3} (X_b) \equiv 0,\  \bmod \ \overline{Y}_{n+3}, \ p+3\leq b \leq n+1.
\end{equation}
This eventually leads to
\begin{equation}
\label{eq:4.4}
\overline{\omega}_{n+3}^1= 0.
\end{equation}
One can check that this change of frame does not affect condition (\ref{eq:4.1}).  We now drop the bars and call this last frame $Y_A$ with Maurer-Cartan forms $\omega_A^B$ satisfying,
\begin{eqnarray}
\label{eq:4.5}
\omega_1^a &=&0, \quad 3 \leq a \leq p+2, \nonumber \\
\omega_{n+3}^b &=&0, \quad p+3 \leq b \leq n+1,\\ 
\omega_1^{n+3} &=& 0, \quad {\omega}_{n+3}^1 = 0. \nonumber
\end{eqnarray}
Furthermore, the following forms are easily shown to be a basis for the cotangent space at each point of $U$,
\begin{equation}
\label{eq:4.6}
\{\omega_{n+3}^3,\ldots,\omega_{n+3}^{p+2},\omega_1^{p+3},\ldots,\omega_1^{n+1} \}.
\end{equation}

We begin by taking the exterior derivative of the equations $\omega_1^a = 0$ and $\omega_{n+3}^b = 0$ in (\ref{eq:4.5}).  Using  (\ref{eq:3.6}),  (\ref{eq:3.7}) and  (\ref{eq:4.5}), we obtain
\begin{equation}
\label{eq:4.7}
0 = \omega_1^{p+3} \wedge \omega_{p+3}^a + \cdots + \omega_1^{n+1} \wedge \omega_{n+1}^a,\quad 3 \leq a \leq p+2,
\end{equation}
\begin{equation}
\label{eq:4.8}
0 = \omega_{n+3}^3 \wedge \omega_3^b + \cdots + \omega_{n+3}^{p+2} \wedge \omega_{p+2}^b,\quad p+3 \leq b \leq n+1.
\end{equation}
We now show that (\ref{eq:4.7}) and (\ref{eq:4.8}) imply that
\begin{equation}
\label{eq:4.9}
\omega_b^a = 0, \quad 3 \leq a \leq p+2, \quad p+3 \leq b \leq n+1.
\end{equation}
To see this, note that since $\omega_b^a = - \omega_a^b$, each of the terms $\omega_b^a$ occurs in exactly one of the equations (\ref{eq:4.7}) and in exactly one of the equations 
(\ref{eq:4.8}).  Equation (\ref{eq:4.7}) involves the basis forms 
$\omega_1^{p+3},\ldots,\omega_1^{n+1}$, while equation (\ref{eq:4.8}) involves the
basis forms $\omega_{n+3}^3,\ldots,\omega_{n+3}^{p+2}$.  We now show how to handle the form 
$\omega_{p+3}^3$; the others are treated in a similar fashion.  The equations from (\ref{eq:4.7}) and (\ref{eq:4.8}), respectively, involving $\omega_{p+3}^3 = - \omega_3^{p+3}$ are
\begin{equation}
\label{eq:4.10}
0 = \omega_1^{p+3} \wedge \omega_{p+3}^3 + \omega_1^{p+4} \wedge \omega_{p+4}^3 + \cdots +  \omega_1^{n+1} \wedge \omega_{n+1}^3,
\end{equation}
\begin{equation}
\label{eq:4.11}
0 = \omega_{n+3}^3 \wedge \omega_3^{p+3} + \omega_{n+3}^4 \wedge \omega_4^{p+3} + \cdots + \omega_{n+3}^{p+2} \wedge \omega_{p+2}^{p+3}.
\end{equation}
We take the wedge product of (\ref{eq:4.10}) with 
$\omega_1^{p+4} \wedge \cdots \wedge \omega_1^{n+1}$ and get
\begin{displaymath}
0 = \omega_{p+3}^3 \wedge (\omega_1^{p+3} \wedge \cdots \wedge \omega_1^{n+1}),
\end{displaymath}
which implies that $\omega_{p+3}$ is in the span of $\omega_1^{p+3},\ldots,\omega_1^{n+1}$.  On the other hand,
taking the wedge product of (\ref{eq:4.11}) with 
$\omega_{n+3}^4 \wedge \cdots \wedge \omega_{n+3}^{p+2}$ yields
\begin{displaymath}
0 = \omega_{p+3}^3 \wedge (\omega_{n+3}^3 \wedge \cdots \wedge \omega_{n+3}^{p+2}),
\end{displaymath}
and thus that $\omega_{p+3}^3$ is in the span of $\omega_{n+3}^3,\ldots,\omega_{n+3}^{p+2}$.  We conclude that $\omega_{p+3}^3 = 0$, as desired.

We next differentiate $\omega_1^{n+3} = 0$ and use (\ref{eq:3.6}), (\ref{eq:3.9}) and (\ref{eq:4.5})
to obtain
\begin{equation}
\label{eq:4.12}
0 = d\omega_1^{n+3} = \omega_1^{p+3} \wedge \omega_{p+3}^{n+3} + \cdots + \omega_1^{n+1} \wedge \omega_{n+1}^{n+3}. 
\end{equation}
This implies that
\begin{equation}
\label{eq:4.13}
\omega_b^{n+3} \in  {\mbox {\rm Span }}\{\omega_1^{p+3},\ldots,\omega_1^{n+1} \}, \quad p+3 \leq b \leq n+1.
\end{equation}
Similarly, differentiation of the equation $\omega_{n+3}^1 = 0$ yields
\begin{equation}
\label{eq:4.14}
0 = \omega_{n+3}^3 \wedge \omega_3^1 + \cdots + \omega_{n+3}^{p+2} \wedge \omega_{p+2}^1, 
\end{equation}
which implies that
\begin{equation}
\label{eq:4.15}
\omega_1^a \in  {\mbox {\rm Span }}\{\omega_{n+3}^3,\ldots,\omega_{n+3}^{p+2} \}, \quad 3 \leq a \leq p+2.
\end{equation}

We next differentiate equation (\ref{eq:4.9}).  Using the skew-symmetry relations (\ref{eq:3.6}) and
equations (\ref{eq:4.5}) and (\ref{eq:4.9}), we see that all terms drop out except the following,
\begin{eqnarray}
0 &=& d \omega_b^a = \omega_b^2 \wedge \omega_2^a + \omega_b^{n+3} \wedge \omega_{n+3}^a \nonumber \\
&=& - (\omega_a^1 \wedge \omega_1^b) + \omega_b^{n+3} \wedge \omega_{n+3}^a. \nonumber 
\end{eqnarray}
Thus, 
\begin{equation}
\label{eq:4.16}
\omega_a^1 \wedge \omega_1^b = \omega_b^{n+3} \wedge \omega_{n+3}^a, \quad 3 \leq a \leq p+2,\quad  p+3 \leq b \leq n+1.
\end{equation}
We now show that equation (\ref{eq:4.16}) implies that there is some function $\alpha$ on $U$ such that
\begin{eqnarray}
\label{eq:4.17}
\omega_a^1 &=& \alpha \  \omega_{n+3}^a, \quad 3 \leq a \leq p+2,\\
\omega_b^{n+3} &=& - \alpha \ \omega_1^b, \quad  p+3 \leq b \leq n+1. \nonumber 
\end{eqnarray}
To see this, note that for any $a$, $3 \leq a \leq p+2$, (\ref{eq:4.15}) gives 
\begin{equation}
\label{eq:4.18}
\omega_a^1 = c_3 \omega_{n+3}^3 + \cdots + c_{p+2} \omega_{n+3}^{p+2}
\  {\rm for} \ {\rm some}  \  c_3, \ldots, c_{p+2},  
\end{equation}
while for any $b$, $p+3 \leq b \leq n+1$, (\ref{eq:4.13}) gives
\begin{equation}
\label{eq:4.19}
\omega_b^{n+3} = d_{p+3} \omega_1^{p+3} + \cdots + d_{n+1} \omega_1^{n+1}
\  {\rm for} \ {\rm some} \  d_{p+3}, \ldots, d_{n+1}.
\end{equation}
Thus,
\begin{equation}
\label{eq:4.20}
\omega_a^1 \wedge \omega_1^b = c_3 \omega_{n+3}^3 \wedge \omega_1^b + \cdots + c_a \omega_{n+3}^a \wedge \omega_1^b
+ \cdots + c_{p+2} \omega_{n+3}^{p+2} \wedge \omega_1^b,
\end{equation}
\begin{equation}
\label{eq:4.21}
\omega_b^{n+3} \wedge \omega_{n+3}^a = d_{p+3} \omega_1^{p+3} \wedge \omega_{n+3}^a + \cdots + d_b \omega_1^b \wedge \omega_{n+3}^a + \cdots + d_{n+1} \omega_1^{n+1} \wedge \omega_{n+3}^a.
\end{equation}
From (\ref{eq:4.16}), we know that the right-hand sides of (\ref{eq:4.20}) and (\ref{eq:4.21})
are equal.  But these expressions contain no common terms from the basis of
2-forms except those involving $\omega_{n+3}^a \wedge \omega_1^b$.  Thence, all of the coefficients except $c_a$ and $d_b$ are zero, and we have
\begin{displaymath}
c_a \omega_{n+3}^a \wedge \omega_1^b = d_b \omega_1^b \wedge \omega_{n+3}^a = (-d_b) \omega_{n+3}^a \wedge \omega_1^b. 
\end{displaymath}
Thus $c_a = - d_b$.  If we set $\alpha_a = c_a$ and $\mu_b = d_b$,  we have shown that
(\ref{eq:4.18}) and (\ref{eq:4.19}) reduce to 
\begin{displaymath}
\omega_a^1 = \alpha_a \omega_{n+3}^a, \quad \omega_b^{n+3} = \mu_b \omega_1^b, 
\  {\rm with} \  \mu_b = - \alpha_a.
\end{displaymath}
This procedure works for any choice of $a$ and $b$ in the appropriate ranges.
By holding $a$ fixed and varying $b$, we see that all of the quantities $\mu_b$ are equal to each other and to $- \alpha_a$.  Similarly, all the quantities $\alpha_a$ are the same, and (\ref{eq:4.17}) holds.  We now consider the expression (\ref{eq:3.5}) for $dY_a, 3 \leq a \leq p+2$. We omit the
terms that vanish because of (\ref{eq:3.6}), (\ref{eq:4.5}) or (\ref{eq:4.9}),
\begin{displaymath}
dY_a = \omega_a^1 Y_1 + \omega_a^3 Y_3 + \cdots + \omega_a^{p+2} Y_{p+2} + \omega_a^{n+2} Y_{n+2} + \omega_a^{n+3} Y_{n+3}.
\end{displaymath}
Using (\ref{eq:4.17}) and the fact from (\ref{eq:3.6}) that $\omega_a^{n+2} = - \omega_{n+3}^a$, this becomes
\begin{equation}
\label{eq:4.22}
dY_a = \omega_{n+3}^a (\alpha Y_1 - Y_{n+2}) + \omega_a^3 Y_3 + \cdots + \omega_a^{p+2} Y_{p+2} + \omega_a^{n+3} Y_{n+3}.
\end{equation}
Similarly, for $p+3 \leq b \leq n+1$, we find
\begin{equation}
\label{eq:4.23}
dY_b = \omega_b^1 Y_1 + \omega_b^2 (Y_2 + \alpha Y_{n+3}) + \omega_b^{p+3} Y_{p+3} + \cdots + \omega_b^{n+1} Y_{n+1}.
\end{equation}
We make the change of frame,
\begin{eqnarray}
\label{eq:4.24}
Y_2^* &=& Y_2 + \alpha Y_{n+3}, \quad Y_{n+2}^* = Y_{n+2} - \alpha Y_1, \\
Y_\beta^* &=& Y_\beta, \quad  \beta \neq 2,\ n+2. \nonumber 
\end{eqnarray}
We now drop the asterisks but use the new frame.  From (\ref{eq:4.22}) and (\ref{eq:4.23}),
we see that in this new frame, we have
\begin{equation}
\label{eq:4.25}
dY_a = \omega_{n+3}^a (- Y_{n+2}) + \omega_a^3 Y_3 + \cdots + \omega_a^{p+2} Y_{p+2} + \omega_a^{n+3} Y_{n+3},
\end{equation}
\begin{equation}
\label{eq:4.26}
dY_b = \omega_b^1 Y_1 + \omega_b^2 Y_2 + \omega_b^{p+3} Y_{p+3} + \cdots + \omega_b^{n+1} Y_{n+1}.
\end{equation}
That is, in the new frame, we have
\begin{equation}
\label{eq:4.27}
\omega_a^1 = 0, \quad 3 \leq a \leq p+2,
\end{equation}
\begin{equation}
\label{eq:4.28}
\omega_b^{n+3} = 0, \quad p+3 \leq b \leq n+1.
\end{equation}
Our goal now is to show that the two spaces
\begin{equation}
\label{eq:4.29}
E = {\mbox {\rm Span }}\{Y_1, Y_2, Y_{p+3}, \ldots, Y_{n+1} \}
\end{equation}
and its orthogonal complement,
\begin{equation}
\label{eq:4.30}
E^{\perp} = {\mbox {\rm Span }}\{Y_3,\ldots,Y_{p+2}, Y_{n+2}, Y_{n+3} \}
\end{equation}
are invariant under $d$.

Concerning $E$, we have that $dY_b \in E$ for $p+3 \leq b \leq n+1$ by (\ref{eq:4.26}).  Furthermore, (\ref{eq:3.6}), (\ref{eq:3.9}) and (\ref{eq:4.5}) imply that 
\begin{displaymath}
dY_1 = \omega_1^1 Y_1 + \omega_1^{p+3} Y_{p+3} + \cdots + \omega_1^{n+1} Y_{n+1},
\end{displaymath}
which is in $E$. Thus, it only remains to show that $dY_2$ is in $E$.  To do this, we differentiate
(\ref{eq:4.27}).  As before, we omit terms which are zero because of
(\ref{eq:3.6}), (\ref{eq:4.5}), (\ref{eq:4.9}) and (\ref{eq:4.27}). We see that formula (\ref{eq:3.7})
for $d\omega_a^1$ reduces to
\begin{equation}
\label{eq:4.31}
0 = d \omega_a^1 = \omega_a^{n+2} \wedge \omega_{n+2}^1 = - \omega_{n+3}^a \wedge \omega_{n+2}^1 = \omega_{n+3}^a \wedge \omega_2^{n+3},\  3 \leq a \leq p+2.
\end{equation}
Similarly, by differentiating equation (\ref{eq:4.28}), we find that
\begin{displaymath}
0 = d \omega_b^{n+3} = \omega_b^2 \wedge \omega_2^{n+3} = - \omega_1^b \wedge \omega_2^{n+3}, \quad p+3 \leq b \leq n+1.
\end{displaymath}
From this and (\ref{eq:4.31}), we see that the wedge product of $\omega_2^{n+3}$ with every form in the basis (\ref{eq:4.6}) is zero, and hence $\omega_2^{n+3} = 0$.  Using this and the fact that $\omega_2^{n+2} = - \omega_{n+3}^1 = 0$, and that by (\ref{eq:3.6}) and (\ref{eq:4.27}), 
\begin{displaymath}
\omega_2^a = - \omega_a^1 = 0, \quad 3 \leq a \leq p+2,
\end{displaymath}
we have
\begin{displaymath}
dY_2 = \omega_2^2 Y_2 + \omega_2^{p+3} Y_{p+3} + \cdots + \omega_2^{n+1} Y_{n+1},
\end{displaymath}
which is in $E$.  So $E$ is invariant under $d$ and is thus a fixed
subspace of $P^{n+2}$, independent of the choice of point of $U$.  Obviously then,
the space $E^{\perp}$ in (\ref{eq:4.30}) is also a fixed subspace of $P^{n+2}$.

Note that $E$ has signature $(q+1,1)$ as a vector subspace of ${\bf R}^{n+3}_2$, and $E^{\perp}$ has signature $(p+1,1)$.  Take an orthornormal basis $\{e_1,\ldots,e_{n+3}\}$ of ${\bf R}^{n+3}_2$ with $e_1$ and $e_{n+3}$ timelike and
\begin{displaymath}
E = {\mbox {\rm Span }}\{e_1,\ldots,e_{q+2}\}, \quad E^{\perp} = 
{\mbox {\rm Span }}\{e_{q+3},\ldots,e_{n+3} \}.
\end{displaymath}
Then $E \cap Q^{n+1}$ is given in homogeneous coordinates $( x_1,\ldots,x_{n+3})$ with respect to this basis by
\begin{equation}
\label{eq:4.32}
x_1^2 = x_2^2 + \cdots + x_{q+2}^2, \quad  x_{q+3} = \cdots = x_{n+3} = 0.
\end{equation}
This quadric is diffeomorphic to the unit sphere $S^q$ in the span $E^{q+1}$ of the spacelike vectors  $e_2,\ldots,e_{q+2}$ with the diffeomorphism $\phi: S^q \rightarrow E \cap Q^{n+1}$ being given by
\begin{equation}
\label{eq:4.33}
\phi (u) = [e_1 + u], \quad u \in S^q.
\end{equation}
Similarly, $E^{\perp} \cap Q^{n+1}$ is the quadric given in homogeneous coordinates by
\begin{equation}
\label{eq:4.34}
x_{n+3}^2 = x_{q+3}^2 + \cdots + x_{n+2}^2, \quad  x_1 = \cdots = x_{q+2} = 0.
\end{equation}
$E^{\perp} \cap Q^{n+1}$ is diffeomorphic to the unit sphere $S^p$ in the span
$E^{p+1}$ of $e_{q+3},\ldots,e_{n+2}$ with the diffeomorphism 
$\psi: S^p \rightarrow E^{\perp} \cap Q^{n+1}$ being given by
\begin{equation}
\label{eq:4.35}
\psi (v) = [v+ e_{n+3}], \quad v \in S^p.
\end{equation}
The focal point map $Y_1$ of our Dupin submanifold is constant on the leaves of the principal foliation $T_r$, and so $Y_1$ factors through an immersion of the $q$-dimensional space of leaves $U/T_r$ into the  $q$-sphere given by the quadric (\ref{eq:4.32}).
Hence, the image of $Y_1$ is an open subset of this quadric. 
Similarly, $Y_{n+3}$ factors through an immersion
of the $p$-dimensional space of leaves $U/T_s$ onto an open subset of the $p$-sphere
given by the quadric (\ref{eq:4.34}).
From this it is clear that the unique compact cyclide containing $\lambda:U \rightarrow \Lambda$ as an open submanifold is given by the Dupin submanifold
$\overline{\lambda} : S^q \times S^p \rightarrow \Lambda$ with
\begin{displaymath}
\overline{\lambda}  (u,v) = [k_1 (u,v), k_2 (u,v)], \quad (u,v) \in S^q \times S^p,
\end{displaymath}
where
\begin{displaymath}
k_1 (u,v) = \phi (u), \quad k_2 (u,v) = \psi (v).
\end{displaymath}
Geometrically, the image of $\overline{\lambda} $ consists of all lines joining a point on the quadric in equation (\ref{eq:4.32}) to a point on the quadric (\ref{eq:4.34}). 

Thus any choice of $(q+1)$-plane $E$ in  $P^{n+2}$ with signature $(q+1,1)$ and corresponding orthogonal complement $E^{\perp}$ determines a unique compact cyclide of characteristic $(p,q)$ and vice-versa. The local Lie equivalence of any two cyclides of the same characteristic is then clear.

From the standpoint of Euclidean geometry, if we consider the point spheres to be given by 
$Q^{n+1} \cap e_{n+3}^\perp$, as in Section 2, then the Legendre submanifold 
$\overline{\lambda}$ above is induced in the usual way from the unit normal bundle
$B^{n-1} = S^q \times S^p$ of the standard embedding of $S^q$ as a great $q$-sphere 
$E^{q+1} \cap S^n$, where $E^{q+1}$ is the span of $e_2, \ldots, e_{q+2}$ and $S^n$ is the unit sphere in $E^{n+1}$, the span of  $e_2, \ldots, e_{n+2}$.  The spheres $S^q$ and 
$S^p = S^n \cap E^{p+1}$, where $E^{p+1}$ is the span of $e_{q+3}, \ldots, e_{n+2}$, are the two focal submanifolds in $S^n$ of a standard product of spheres $S^p \times S^q$ in $S^n$ (see
\cite[pp. 295--296]{CR7} or \cite[pp. 110--111]{CR8}).\\

\section{Dupin Submanifolds for $n = 4$.}
\label{sec:5}

The classification of Dupin submanifolds induced from surfaces in ${\bf R}^3$ follows from the results of the last section.  In his doctoral dissertation \cite{P1} (later published as \cite{P3}), U. Pinkall obtained a local classification up to Lie equivalence of all Dupin submanifolds induced from hypersurfaces in ${\bf R}^4$.  As we shall see, this is a far more complicated calculation than that of the previous section, and as yet, no one has obtained a similar classification of Dupin hypersurfaces in ${\bf R}^n$ for $n \geq 5$.  In this section, we will prove Pinkall's theorem using the method of moving frames.\\  

\noindent
{\bf Remark 5.1.} (added by T. Cecil in 2020) See \cite[pp. 301--308]{CR8} for a survey of local classifications of higher dimensional proper Dupin hypersurfaces that have been found since the time that the original version of this paper was written. \\

We follow the notation used in Sections 3 and 4.  We consider a Dupin submanifold
\begin{equation}
\label{eq:5.1}
\lambda:B \rightarrow \Lambda
\end{equation}
where
\begin{equation}
\label{eq:5.2}
\dim B = 3, \quad \dim \Lambda = 7,
\end{equation}
and the image $\lambda(b), b \in B$, is the line $[Y_1, Y_7]$ of the Lie frame
$Y_1,\ldots, Y_7$.  We assume that there are three distinct focal points on each line $\lambda(b)$.

By (\ref{eq:3.39}), we can choose the frame so that
\begin{eqnarray}
\label{eq:5.3}
&\omega_1^3 = \omega_1^7 = 0,& \\
&\omega_7^4 = \omega_7^1 = 0.&\nonumber
\end{eqnarray}
By making a change of frame of the form
\begin{eqnarray}
\label{eq:5.4}
&Y_1^* = \alpha Y_1, \quad Y_2^* = (1/\alpha) Y_2, & \\
&Y_7^* = \beta Y_7, \quad Y_6^* = (1/\beta) Y_6,& \nonumber
\end{eqnarray}
for suitable smooth functions $\alpha$ and $\beta$ on $B$, we can arrange that $Y_1 + Y_7$ represents the third focal point for each point of $B$.  Then by using the fact that $B$ is Dupin, we can use the method employed at the end of Section 3 to make a change of frame leading to the following equations similar to (\ref{eq:3.39}) (and to (\ref{eq:5.3})),
\begin{eqnarray}
\label{eq:5.5}
&\omega_1^5 + \omega_7^5 = 0,& \\ 
&\omega_1^1 - \omega_7^7 = 0.& \nonumber
\end{eqnarray}
This completely fixes the $Y_i, i = 3,4,5$, and $Y_1, Y_7$ are determined up to a transformation of the form
\begin{equation}
\label{eq:5.6}
Y_1^* = \tau Y_1, \quad Y_7^* = \tau Y_7.
\end{equation}

Each of the three focal point maps $Y_1$, $Y_7$, $Y_1 + Y_7$ is constant along the leaves of its corresponding principal foliation.  Thus, each focal point map factors through an immersion of the corresponding 2-dimensional space of leaves of its principal foliation into $Q^5$.  (See Section 4 of Chapter 2 of the book \cite{CR7} or \cite[pp. 18--32]{CR8} for more detail on this point.)  In terms of moving frames, this implies that the forms $\omega_1^4, \omega_1^5, \omega_7^3$ are linearly independent on $B$, i.e.,
\begin{equation}
\label{eq:5.7}
\omega_1^4 \wedge \omega_1^5 \wedge \omega_7^3 \neq 0.
\end{equation}

This can also be seen by expressing the forms above in terms of a Lie frame 
$\{\overline{Y}_1,\ldots,\overline{Y}_{n+3}\}$, where $\overline{Y}_1$
satisfies the regularity condition (\ref{eq:3.8}), and using the fact that each focal point has multiplicity one.  For simplicity, we will also use the notation,
\begin{equation}
\label{eq:5.8}
\theta_1 = \omega_1^4, \quad \theta_2 = \omega_1^5, \quad \theta_3 = \omega_7^3.
\end{equation}

Analytically, the Dupin conditions are three partial differential equations, and we are treating an over-determined system.  The method of moving frames reduces the handling of its integrability conditions to a straightforward
algebraic problem, viz., that of repeated exterior differentiations.

We begin by computing the exterior derivatives of the three equations 
$\omega_1^3 = 0,$ $\omega_7^4 = 0,$ $\omega_1^5 + \omega_7^5 = 0$.
Using the skew-symmetry relations (\ref{eq:3.6}), as well as (\ref{eq:5.3}) and (\ref{eq:5.5}), the exterior derivatives of these three equations yield the system
\begin{eqnarray}
\label{eq:4.6.7}
0 & = & \omega_1^4 \wedge \omega_3^4 + \omega_1^5 \wedge \omega_3^5,\nonumber \\
0 & = & \omega_1^5 \wedge \omega_4^5 + \omega_7^3 \wedge \omega_3^4, \nonumber \\
0 & = & \omega_1^4 \wedge \omega_4^5 + \omega_7^3 \wedge \omega_3^5.\nonumber
\end{eqnarray}
If we take the wedge product of the first of these with $\omega_1^4$, we conclude that
$ \omega_3^5$ is in the span of $\omega_1^4$ and $\omega_1^5$.  On the other hand, taking the wedge product of of the third equation
with $\omega_1^4$ yields that $\omega_3^5$ is in the span of $\omega_1^4$ and $\omega_7^3$.  Consequently, $\omega_3^5 = \rho \omega_1^4$, for some smooth function $\rho$ on $B$.  Similarly, one can show that there exist smooth functions $\sigma$
and $\tau$ such that $\omega_3^4 = \sigma \omega_1^5$ and $\omega_4^5 = \tau \omega_7^3$.  Then, if we substitute these into the three equations above, we get that $\rho = \sigma = \tau$, and hence we have,
\begin{equation}
\label{eq:5.9}
\omega_3^5 = \rho \omega_1^4, \quad \omega_3^4 = \rho \omega_1^5, \quad \omega_4^5 = \rho \omega_7^3.
\end{equation}
Next we differentiate the equations $\omega_1^7 = 0,$ $\omega_7^1 = 0,$ 
$\omega_1^1 - \omega_7^7 = 0.$  
As above, use of the skew-symmetry relations (\ref{eq:3.6}) and the equations
(\ref{eq:5.3}) and (\ref{eq:5.5}) yields the existence
of smooth functions $a$, $b$, $c$, $p$, $q$, $r$, $s$, $t$, $u$ on $B$ such that the following relations hold:
\begin{eqnarray}
\label{eq:5.10}
\omega_4^7 & = & - \omega_6^4 = a \omega_1^4 + b \omega_1^5, \\
\omega_5^7 & = & - \omega_6^5 = b \omega_1^4 + c \omega_1^5;\nonumber
\end{eqnarray}
\begin{eqnarray}
\label{eq:5.11}
\omega_3^1 & = & - \omega_2^3 = p \omega_7^3 - q \omega_1^5, \\
\omega_5^1 & = & - \omega_2^5 = q \omega_7^3 - r \omega_1^5;\nonumber
\end{eqnarray}
\begin{eqnarray}
\label{eq:5.12}
\omega_4^1 & = & - \omega_2^4 = b \omega_1^5 + s \omega_1^4 + t \omega_7^3, \\
\omega_6^3 & = & - \omega_3^7 = q \omega_1^5 + t \omega_1^4 + u \omega_7^3.\nonumber
\end{eqnarray}

We next see what can be deduced from taking the exterior derivatives of the equations
(\ref{eq:5.9})--(\ref{eq:5.12}).  First, we take the exterior derivatives of the three basis forms
$\omega_1^4, \omega_1^5, \omega_7^3$.  For example, using the relations that we have derived so far, we have from the Maurer-Cartan equation (\ref{eq:3.7}), 
\begin{displaymath}
d\omega_1^4 = \omega_1^1 \wedge \omega_1^4 + \omega_1^5 \wedge \omega_5^4 = \omega_1^1 \wedge \omega_1^4 
- \rho \omega_1^5 \wedge \omega_7^3.
\end{displaymath}
We obtain similar equations for $d\omega_1^5$ and $d\omega_7^3$. When we use the forms
$\theta_1$, $\theta_2$, $\theta_3$ defined in (\ref{eq:5.8}) for $\omega_1^4, \omega_1^5, \omega_7^3$,
we have
\begin{eqnarray}
\label{eq:5.13}
d\theta_1 & = & \omega_1^1 \wedge \theta_1 - \rho \ \theta_2 \wedge \theta_3, \nonumber\\
d\theta_2 & = & \omega_1^1 \wedge \theta_2 - \rho \ \theta_3 \wedge \theta_1, \\
d\theta_3 & = & \omega_1^1 \wedge \theta_3 - \rho \ \theta_1 \wedge \theta_2. \nonumber
\end{eqnarray}
We next differentiate (\ref{eq:5.9}).  We have $\omega_3^4 = \rho \omega_1^5$.  On the one hand,
\begin{displaymath}
d\omega_3^4 = \rho \ d\omega_1^5 + d\rho \wedge \omega_1^5.
\end{displaymath}
Using the second equation in (\ref{eq:5.13}) with $\omega_1^5 = \theta_2$, this becomes
\begin{displaymath}
d\omega_3^4 = \rho \omega_1^1 \wedge \omega_1^5 - \rho^2 \omega_7^3 \wedge \omega_1^4 + d\rho \wedge \omega_1^5.
\end{displaymath}
On the other hand, we can compute $d\omega_3^4$ from the Maurer-Cartan equation 
(\ref{eq:3.7}) and use the relationships that we have derived to find
\begin{displaymath}
d\omega_3^4 = (-p - \rho^2 -a) (\omega_1^4 \wedge \omega_7^3) - q \omega_1^5 \wedge \omega_1^4 + b \omega_7^3 \wedge \omega_1^5.
\end{displaymath}
If we equate these two expressions for $d\omega_3^4$, we get
\begin{equation}
\label{eq:5.14}
(-p -a -2\rho^2)\  \omega_1^4 \wedge \omega_7^3 = (d\rho + \rho \omega_1^1 - q \omega_1^4 - b \omega_7^3) \wedge \omega_1^5.
\end{equation}
Because of the independence of $\omega_1^4, \omega_1^5$ and 
$\omega_7^3$, both sides of the equation above must vanish.  Thus, we conclude that
\begin{equation}
\label{eq:5.15}
2 \rho^2 = -a-p,
\end{equation}
and that $d\rho + \rho \omega_1^1 - q \omega_1^4 - b \omega_7^3$ is a multiple of $\omega_1^5$.  Similarly, differentiation of the equation $\omega_4^5 = \rho \omega_7^3$ yields the following analogue of (\ref{eq:5.14}),
\begin{equation}
\label{eq:5.16}
(s-a-r+2\rho^2)\  \omega_1^4 \wedge \omega_1^5 = (d\rho + \rho \omega_1^1 +t \omega_1^5 - q \omega_1^4) \wedge \omega_7^3,
\end{equation}
and differentiation of $\omega_3^5 = \rho \omega_1^4$ yields
\begin{equation}
\label{eq:5.17}
(c+p+u-2\rho^2)\  \omega_1^5 \wedge \omega_7^3 = (-d\rho - \rho \omega_1^1 - t \omega_1^5 +b \omega_7^3) \wedge \omega_1^4.
\end{equation}
In each of the equations (\ref{eq:5.14}), (\ref{eq:5.16}), (\ref{eq:5.17}), both sides of the equation must vanish.
From the vanishing of the left sides of the equations, we get the fundamental relationship,
\begin{equation}
\label{eq:5.18}
2\rho^2= -a-p = a+r-s = c+p+u.
\end{equation}
Furthermore, from the vanishing of the right-hand sides of the three equations 
(\ref{eq:5.14}), (\ref{eq:5.16}) and (\ref{eq:5.17}),
we can determine after some algebra that
\begin{equation}
\label{eq:5.19}
d\rho + \rho \omega_1^1 = q \omega_1^4 - t \omega_1^5 + b \omega_7^3.
\end{equation}
The last equation shows the importance of $\rho$. 
Following the notation introduced in equation (\ref{eq:5.8}), we write equation (\ref{eq:5.19}) as,
\begin{equation}
\label{eq:5.20}
d\rho + \rho \omega_1^1 =  \rho_1 \theta_1 + \rho_2 \theta_2 + \rho_3 \theta_3,
\end{equation}
where
\begin{equation}
\label{eq:5.21}
\rho_1 = q, \quad \rho_2 = -t, \quad \rho_3 = b,
\end{equation}
are the ``covariant derivatives'' of $\rho$.

Using the Maurer-Cartan equations, we can compute
\begin{eqnarray}
d\omega_1^1 & = & \omega_1^4 \wedge \omega_4^1 + \omega_1^5 \wedge \omega_5^1\nonumber \\
& = & \omega_1^4 \wedge (b \omega_1^5 + t \omega_7^3) + \omega_1^5 \wedge (q \omega_7^3 - r \omega_1^5) \nonumber\\
& = & b  \omega_1^4 \wedge \omega_1^5 + q \omega_1^5 \wedge \omega_7^3 - t \omega_7^3 \wedge \omega_1^4.\nonumber
\end{eqnarray}
Using (\ref{eq:5.8}) and (\ref{eq:5.21}), this can be rewritten as
\begin{equation}
\label{eq:5.22}
d\omega_1^1 = \rho_3 \ \theta_1 \wedge \theta_2 + \rho_1 \ \theta_2 \wedge \theta_3 + \rho_2 \ \theta_3 \wedge \theta_1.
\end{equation}
The trick now is to express everything in terms of $\rho$ and its successive covariant
derivatives. 

We first derive a general form for these covariant derivatives.  Suppose that $\sigma$ is a smooth function which satisfies a relation of the form
\begin{equation}
\label{eq:5.23}
d\sigma + m \sigma \omega_1^1 = \sigma_1 \ \theta_1 + \sigma_2 \  \theta_2 + \sigma_3 \  \theta_3,
\end{equation}
for some integer $m$. (Note that (\ref{eq:5.19}) is such a relationship for the function $\rho$
with $m = 1$.)  By taking the exterior derivative of (\ref{eq:5.23}) and using
(\ref{eq:5.13}) and (\ref{eq:5.22}) to express both sides in terms of the standard basis
of 2-forms $\theta_1 \wedge \theta_2$, $\theta_2 \wedge \theta_3$ and $\theta_3 \wedge \theta_1$, one finds
that the functions $\sigma_1, \sigma_2, \sigma_3$ satisfy equations of the form
\begin{equation}
\label{eq:5.24}
d\sigma_{\alpha} + (m+1) \sigma_{\alpha} \omega_1^1 = \sigma_{\alpha 1} \ \theta_1 
+ \sigma_{\alpha 2} \  \theta_2 + \sigma_{\alpha 3} \  \theta_3, \quad \alpha = 1,2,3,
\end{equation}
where the coefficient functions $\sigma_{\alpha \beta}$ satisfy the commutation relations,\index{commutation relations on covariant derivatives}
\begin{eqnarray}
\label{eq:5.25}
\sigma_{12} - \sigma_{21} & = & - m \sigma \rho_3 - \rho \sigma_3 , \nonumber\\
\sigma_{23} - \sigma_{32} & = & - m \sigma \rho_1 - \rho \sigma_1 , \\
\sigma_{31} - \sigma_{13} & = & - m \sigma \rho_2 - \rho \sigma_2 . \nonumber
\end{eqnarray}
In particular, from equation (\ref{eq:5.20}), we have the following commutation relations on $\rho_1, \rho_2, \rho_3$:
\begin{eqnarray}
\label{eq:5.26}
\rho_{12} - \rho_{21} & = & - 2\rho \rho_3 \nonumber\\
\rho_{23} - \rho_{32} & = & - 2\rho \rho_1\\
\rho_{31} - \rho_{13} & = & - 2\rho \rho_2. \nonumber
\end{eqnarray}
We next take the exterior derivative of equations (\ref{eq:5.10})--(\ref{eq:5.12}).  We first differentiate the equation
\begin{equation}
\label{eq:5.27}
\omega_4^7 = a \omega_1^4 + b \omega_1^5.
\end{equation}
On the one hand, from the Maurer-Cartan equation (\ref{eq:3.7}) for $d\omega_4^7$, we have
(by not writing those terms which have already been shown to vanish),
\begin{eqnarray}
\label{eq:5.28}
d\omega_4^7 & = & \omega_4^2 \wedge \omega_2^7 + \omega_4^3 \wedge \omega_3^7 + \omega_4^5 \wedge \omega_5^7 + 
\omega_4^7 \wedge \omega_7^7 \\
& = &  - \omega_1^4 \wedge \omega_2^7 + (-\rho \omega_1^5) \wedge (-q \omega_1^5 -t \omega_1^4 - u \omega_7^3) \nonumber\\
& + & \rho \omega_7^3 \wedge (b \omega_1^4 + c \omega_1^5) +(a \omega_1^4 + b \omega_1^5) \wedge \omega_1^1.\nonumber
\end{eqnarray}
On the other hand, differentiation of the right-hand side of equation (\ref{eq:5.27}) yields
\begin{eqnarray}
\label{eq:5.29}
d\omega_4^7 & = & da \wedge \omega_1^4 + a\  d\omega_1^4  + db \wedge \omega_1^5 + b\  d\omega_1^5  \\
& = &  da \wedge \omega_1^4 + a (\omega_1^1 \wedge \omega_1^4 - \rho \omega_1^5 \wedge \omega_7^3) \nonumber\\
& + & db \wedge \omega_1^5 + b (\omega_1^1 \wedge \omega_1^5 - \rho \omega_1^4 \wedge \omega_7^3).  \nonumber
\end{eqnarray}
Equating (\ref{eq:5.28}) and (\ref{eq:5.29}) yields
\begin{eqnarray}
\label{eq:5.30}
(da &+& 2a \omega_1^1 - 2b \rho \omega_7^3 - \omega_2^7) \wedge \omega_1^4 \\
&+&  (db + 2b \omega_1^1 +(a+u-c)\rho \omega_7^3) \wedge \omega_1^5 
+ \rho t \omega_1^4 \wedge \omega_1^5 = 0.  \nonumber
\end{eqnarray}
Since $b = \rho_3$, it follows from (\ref{eq:5.19}) and (\ref{eq:5.24}) that
\begin{equation}
\label{eq:5.31}
db + 2b \omega_1^1 = d\rho_3 + 2 \rho_3 \omega_1^1 = \rho_{31}\ \theta_1 + \rho_{32}\ \theta_2 + \rho_{33} \ \theta_3.
\end{equation}
By examining the coefficient of $\omega_1^5 \wedge \omega_7^3 = \theta_2 \wedge \theta_3$ in equation  (\ref{eq:5.30}) and using (\ref{eq:5.31}), we get that
\begin{equation}
\label{eq:5.32}
\rho_{33} = \rho (c-a-u).
\end{equation}
Furthermore, the remaining terms in (\ref{eq:5.30}) are
\begin{eqnarray}
\label{eq:5.33}
&  & (da + 2a \omega_1^1 - \omega_2^7 - 2\rho b \omega_7^3 - (\rho t + \rho_{31}) \omega_1^5) \wedge \omega_1^4  \\
& + & {\mbox{\rm terms involving }} \omega_1^5 \ {\mbox{\rm and }} \omega_7^3 \ {\mbox{\rm only.}}\nonumber
\end{eqnarray}
Thus, the coefficient in parentheses must be a multiple of $\omega_1^4$, call it 
$\bar{a} \omega_1^4$.  
We can write this using (\ref{eq:5.8}) and (\ref{eq:5.21}) as
\begin{equation}
\label{eq:5.34}
da + 2a \omega_1^1 = \omega_2^7 + \bar{a} \theta_1 + (\rho_{31} - \rho \rho_2)\theta_2
+ 2\rho \rho_3 \theta_3.
\end{equation}
In a similar manner, if we differentiate
\begin{displaymath}
\omega_5^7 = b \omega_1^4 + c \omega_1^5,
\end{displaymath}
we obtain,
\begin{equation}
\label{eq:5.35}
dc + 2c \omega_1^1 = \omega_2^7 + (\rho_{32} + \rho \rho_1) \theta_1 + \bar{c} \theta_2
- 2 \rho \rho_3 \theta_3.
\end{equation}
Thus, from the two equations in (\ref{eq:5.10}), we have obtained equations (\ref{eq:5.32}), 
(\ref{eq:5.34}) and (\ref{eq:5.35}).  In completely analogous fashion, we can differentiate the two equations in (\ref{eq:5.11}) to obtain
\begin{equation}
\label{eq:5.36}
\rho_{11} = \rho (s+r-p),
\end{equation}
\begin{equation}
\label{eq:5.37}
dp + 2p \omega_1^1 = - \omega_2^7 + 2 \rho \rho_1 \theta_1 + (- \rho_{13} - \rho \rho_2) \theta_2
+ \bar{p} \theta_3,
\end{equation}
\begin{equation}
\label{eq:5.38}
dr + 2r \omega_1^1 = - \omega_2^7 - 2 \rho \rho_1 \theta_1 + \bar{r} \theta_2 + 
(- \rho_{12} + \rho \rho_3) \theta_3,
\end{equation}
while differentiation of (\ref{eq:5.12}) yields
\begin{equation}
\label{eq:5.39}
\rho_{22} + \rho_{33} = \rho (p-r-s),
\end{equation}
\begin{equation}
\label{eq:5.40}
ds + 2s \omega_1^1 = \bar{s} \theta_1 + (\rho_{31} + \rho \rho_2) \theta_2 + (-\rho_{21} + \rho \rho_3) 
\theta_3,
\end{equation}
\begin{equation}
\label{eq:5.41}
du + 2u \omega_1^1 = (- \rho_{23} - \rho \rho_1) \theta_1 + (\rho_{13} - \rho \rho_2) \theta_2 + \bar{u} \theta_3.
\end{equation}

In these equations, the coefficients $\bar{a}$, $\bar{c}$, $\bar{p}$, $\bar{r}$, $\bar{s}$, $\bar{u}$ remain undetermined.  However, by differentiating equation (\ref{eq:5.18}) and using the appropriate equations among those involving these quantities above, one can show that
\begin{eqnarray}    
\label{eq:5.42}
\bar{a} & = &- 6 \rho \rho_1, \quad \bar{c} = 6 \rho \rho_2,\nonumber \\
\bar{p} & = &- 6 \rho \rho_3, \quad \bar{r} = 6 \rho \rho_2, \\
\bar{s} & = &- 12 \rho \rho_1, \quad \bar{u} = 12 \rho \rho_3.\nonumber 
\end{eqnarray}
From equations (\ref{eq:5.32}), (\ref{eq:5.36}), (\ref{eq:5.39}) and (\ref{eq:5.18}), we easily compute that 
\begin{equation}
\label{eq:5.43}
\rho_{11} + \rho_{22} + \rho_{33} = 0.
\end{equation}
Using (\ref{eq:5.42}), equations (\ref{eq:5.40}) and (\ref{eq:5.41}) can be rewritten as
\begin{equation}
\label{eq:5.44}
ds + 2s \omega_1^1 = - 12 \rho \rho_1 \theta_1 +(\rho_{31} + \rho \rho_2) \theta_2 + (-\rho_{21} + \rho \rho_3) \theta_3,
\end{equation}
\begin{equation}
\label{eq:5.45}
du + 2u \omega_1^1 = (-\rho_{23} - \rho \rho_1) \theta_1 +(\rho_{13} - \rho \rho_2) \theta_2 + 12 \rho \rho_3 \theta_3.
\end{equation}
By taking the exterior derivatives of these two equations and making use of equation (\ref{eq:5.43})
and of the commutation relations (\ref{eq:5.25}) for $\rho$ and its various derivatives,
one can ultimately show after a lengthy calculation that the following fundamental equations hold:
\begin{eqnarray}
\label{eq:5.46}
\rho \rho_{12} & + & \rho_1 \rho_2 + \rho^2 \rho_3 = 0,\nonumber \\
\rho \rho_{21} & + & \rho_1 \rho_2 - \rho^2 \rho_3 = 0,\nonumber \\
\rho \rho_{23} & + & \rho_2 \rho_3 + \rho^2 \rho_1 = 0, \\
\rho \rho_{32} & + & \rho_2 \rho_3 - \rho^2 \rho_1 = 0,\nonumber \\
\rho \rho_{31} & + & \rho_3 \rho_1 + \rho^2 \rho_2 = 0,\nonumber \\
\rho \rho_{13} & + & \rho_3 \rho_1 - \rho^2 \rho_2 = 0.\nonumber 
\end{eqnarray}

We now briefly outline the details of this calculation.  By (\ref{eq:5.44}), we have
\begin{equation}
\label{eq:5.47}
s_1 = -12 \rho \rho_1, \quad s_2 = \rho_{31} + \rho \rho_2, \quad s_3 = \rho \rho_3 - \rho_{21}.
\end{equation}
The commutation relation (\ref{eq:5.25}) for $s$ with $m=2$ gives
\begin{equation}
\label{eq:5.48}
s_{12} - s_{21} = - 2s \rho_3 - \rho s_3 = - 2s \rho_3 - \rho (\rho \rho_3 - \rho_{21}).
\end{equation}
On the other hand, we can directly compute by taking covariant
derivatives of (\ref{eq:5.47}) that
\begin{equation}
\label{eq:5.49}
s_{12} - s_{21} = -12 \rho \rho_{12} - 12 \rho_2 \rho_1 - (\rho_{311} + \rho_1 \rho_2 + \rho \rho_{21}).
\end{equation}
The main problem now is to get $\rho_{311}$ into a usable form. By taking the covariant
derivative of the third equation in (\ref{eq:5.26}), we find
\begin{equation}
\label{eq:5.50}
\rho_{311} - \rho_{131} = -2 \rho_1 \rho_2 - 2 \rho \rho_{21}.
\end{equation}
Then using the commutation relation,
\begin{displaymath}
\rho_{131} = \rho_{113} - 2 \rho_1 \rho_2 - \rho \rho_{12},
\end{displaymath}
we get from (\ref{eq:5.50}) that
\begin{equation}
\label{eq:5.51}
\rho_{311} = \rho_{113} - 4 \rho_1 \rho_2 - \rho \rho_{12} - 2\rho \rho_{21}.
\end{equation}
Taking the covariant derivative of  $\rho_{11} = \rho (s+r-p)$ and substituting the expression obtained for $\rho_{113}$ into (\ref{eq:5.51}), we get
\begin{equation}
\label{eq:5.52}
\rho_{311} = \rho_3 (s+r-p) - 3 \rho \rho_{21} - 2 \rho \rho_{12} + 8 \rho^2 \rho_3 - 4 \rho_1 \rho_2.
\end{equation}
If we substitute (\ref{eq:5.52}) for $\rho_{311}$ in (\ref{eq:5.49}) and then equate the right-hand  sides of (\ref{eq:5.48}) and (\ref{eq:5.49}), we obtain the first equation in (\ref{eq:5.46}).  The cyclic
permutations are obtained in a similar way from $s_{23} - s_{32}$, etc.

Our frame attached to the line $[Y_1,Y_7]$ is still not completely determined, viz., the following change is allowable:
\begin{equation}
\label{5.53}
Y_2^* = \alpha^{-1}  Y_2 + \mu Y_7, \quad Y_6^* =   \alpha^{-1} Y_6 - \mu Y_1.
\end{equation}
The $Y_i$'s, $i=3,4,5$ being completely determined , we have under this change,
\begin{eqnarray}
\omega_1^{4*} & = & \alpha \omega_1^4, \quad  \omega_1^{5*} = \alpha \omega_1^5, \quad  \omega_7^{3*} = \alpha \omega_7^3,
\nonumber \\
\omega_4^{7*} & = & \alpha^{-1} \omega_4^7 + \mu \omega_1^4, \nonumber \\
\omega_3^{1*} & = & \alpha^{-1} \omega_3^1 - \mu \omega_7^3, \nonumber
\end{eqnarray}
which implies that
\begin{displaymath}
a^* = \alpha^{-2} a + \alpha^{-1} \mu, \quad p^* = \alpha^{-2} p - \alpha^{-1} \mu.
\end{displaymath}
We choose $\mu$ to make $a^* = p^*$.  After dropping the asterisks, we have from
(\ref{eq:5.18}) that
\begin{equation}
\label{eq:5.54}
a=p= - \rho^2, \quad r = 3 \rho^2 +s, \quad c = 3 \rho^2 - u.
\end{equation}
Now using the fact that $a=p$, we can subtract (\ref{eq:5.37}) from (\ref{eq:5.34}) and get that
\begin{equation}
\label{eq:5.55}
\omega_2^7 = 4 \rho \rho_1 \theta_1 - ((\rho_{31}+\rho_{13})/2) \theta_2 - 4 \rho \rho_3 \theta_3.
\end{equation}

We are finally in position to proceed toward the main results.
Ultimately, we show that the frame can be chosen so that the function $\rho$ is constant, and the  classification naturally splits into two cases $\rho = 0$ and $\rho \neq 0$.\\

\noindent
{\em The case $\rho \neq 0$.}\\

\noindent
We now assume  that the function $\rho$ is never zero on $B$.  The following lemma is the key in this case.  This is Pinkall's Lemma \cite[p. 108]{P3}, where his function $c$ is the negative of our function $\rho$.  Since $\rho \neq 0$, the fundamental equations (\ref{eq:5.46}) allow one  to express all of the second covariant derivatives $\rho_{\alpha \beta}$ in terms of $\rho$ and its first derivatives $\rho_\alpha$.  This enables us to give a somewhat simpler proof than Pinkall gave for the lemma.\\

\noindent
{\bf Lemma 5.1.} {\it Suppose that $\rho$ never vanishes on  $B$.  Then 
$\rho_1 = \rho_2 = \rho_3 = 0$ at every point of $B$.}\\

\noindent
{\bf Proof:}  First note that if the function $\rho_3$ vanishes identically, then (\ref{eq:5.46}) and the assumption that $\rho \neq 0$ imply that $\rho_1$ and $\rho_2$ also vanish identically.  We now complete the proof by showing that $\rho_3$ must vanish everywhere.  This is accomplished by considering the expression $s_{12} - s_{21}$.
By the commutation relations (\ref{eq:5.25}), we have
\begin{displaymath}
s_{12} - s_{21} = -2s \rho_3 - \rho s_3.
\end{displaymath}
By (\ref{eq:5.46}) and (\ref{eq:5.47}), we see that
\begin{displaymath}
\rho s_3 = \rho^2 \rho_3 - \rho \rho_{21} = \rho_1 \rho_2,
\end{displaymath} 
and so
\begin{equation}
\label{eq:5.56}
s_{12} - s_{21} = -2s \rho_3 - \rho_1 \rho_2.
\end{equation}
On the other hand, we can compute $s_{12}$ directly from the equation
$s_1 = -12 \rho \rho_1$. 
Then using the expression for $\rho_{12}$ obtained from (\ref{eq:5.46}), we get
\begin{eqnarray}
\label{eq:5.57}
s_{12} & = & -12 \rho_2 \rho_1 - 12 \rho \rho_{12} = -12 (\rho_2 \rho_1 +\rho \rho_{12}) \\
& = & -12 (\rho_2 \rho_1 + (-\rho_2 \rho_1 - \rho^2 \rho_3)) = 12 \rho^2 \rho_3. \nonumber
\end{eqnarray}
Next we have from (\ref{eq:5.47}),
\begin{displaymath}
s_2 = \rho_{31} + \rho \rho_2.  
\end{displaymath}
Using (\ref{eq:5.46}), we can write
\begin{displaymath}
\rho_{31} = -\rho_3 \rho_1 \rho^{-1} - \rho \rho_2,
\end{displaymath} 
and thus,
\begin{equation}
\label{eq:5.58}
s_2 = - \rho_3 \rho_1 / \rho.
\end{equation}
Then, we compute
\begin{displaymath}
s_{21} = - (\rho (\rho_3 \rho_{11} + \rho_{31} \rho_1) - \rho_3 \rho_1^2)/ \rho^2.
\end{displaymath} 
Using (\ref{eq:5.36}) for $\rho_{11}$ and (\ref{eq:5.46}) to get $\rho_{31}$, this becomes
\begin{equation}
\label{eq:5.59}
s_{21} = - \rho_3 (s+r-p) + 2 \rho_3 \rho_1^2 \rho^{-2} + \rho_1 \rho_2.
\end{equation}
Now equate the expression in (\ref{eq:5.56}) for $s_{12} - s_{21}$ with that obtained by subtracting
equation (\ref{eq:5.59}) from equation (\ref{eq:5.57}) to get
\begin{displaymath}
- 2s \rho_3 - \rho_1 \rho_2 = 12 \rho^2 \rho_3 + \rho_3 (s+r-p) - 2 \rho_3 \rho_1^2 \rho^{-2} - \rho_1 \rho_2.
\end{displaymath}
This can be rewritten as
\begin{equation}
\label{eq:5.60}
0 = \rho_3 (12 \rho^2 + 3s + r -p - 2 \rho_1^2 \rho^{-2} ).
\end{equation}
Using the expressions in (\ref{eq:5.54}) for $r$ and $p$, we see that
\begin{displaymath}
3s+r-p = 4s + 4 \rho^2,
\end{displaymath}
and so (\ref{eq:5.60}) can be written as
\begin{equation}
\label{eq:5.61}
0 = \rho_3 (16 \rho^2 + 4s - 2 \rho_1^2 \rho^{-2} ).
\end{equation}
Suppose that $\rho_3 \neq 0$ at some point $b \in B$.  Then $\rho_3$ does not vanish on some neighborhood $U$ of $b$.  By (\ref{eq:5.61}), we have
\begin{equation}
\label{eq:5.62}
16 \rho^2 + 4s - 2 \rho_1^2 \rho^{-2}  = 0
\end{equation}
on $U$.  We now take the $\theta_2$-covariant derivative of (\ref{eq:5.62}) and obtain
\begin{equation}
\label{eq:5.63}
32 \rho \rho_2 + 4s_2 - 4 \rho_1 \rho_{12} \rho^{-2} + 4 \rho_1^2 \rho_2 \rho^{-3} =0.
\end{equation}
We now substitute the expression (\ref{eq:5.58}) for $s_2$ and the formula
\begin{displaymath}
\rho_{12} = - \rho_1 \rho_2 \rho^{-1} - \rho \rho_3
\end{displaymath}
obtained from (\ref{eq:5.46}) into (\ref{eq:5.63}).  After some algebra, (\ref{eq:5.63})
reduces to
\begin{displaymath}
\rho_2 (32 \rho^4 + 8 \rho_1^2) = 0.
\end{displaymath}
Since $\rho \neq 0$, this implies that $\rho_2 = 0$ on $U$.  But then the left side of the equation 
(\ref{eq:5.46})
\begin{displaymath}
\rho \rho_{21} + \rho_1 \rho_2 = \rho^2 \rho_3,
\end{displaymath}
must vanish on $U$.  Since $\rho \neq 0$, we conclude that $\rho_3 = 0$ on $U$, a contradiction to our assumption. Hence, $\rho_3$ vanishes identically on $B$, and the lemma is proven.

We now continue with the case $\rho \neq 0$.  According to Lemma 5.1, all the covariant derivatives of $\rho$ are zero, and our formulas simplify greatly.  Equations (\ref{eq:5.32}) and (\ref{eq:5.36}) give
\begin{displaymath}
c-a-u = 0, \quad s+r-p =0.
\end{displaymath}
These combined with (\ref{eq:5.54}) give
\begin{equation}
\label{eq:5.64}
c = r = \rho^2, \quad u = -s = 2\rho^2.
\end{equation}
By (\ref{eq:5.55}) we have $\omega_2^7 = 0$.  So the differentials of the frame vectors can now be written
\begin{eqnarray}
\label{eq:5.65}
dY_1  -  \omega_1^1 Y_1 & = & \omega_1^4 Y_4 + \omega_1^5 Y_5, \nonumber \\
dY_7  -  \omega_1^1 Y_7 & = & \omega_7^3 Y_3 - \omega_1^5 Y_5, \nonumber \\
dY_2  +  \omega_1^1 Y_2 & = & \rho^2 (\omega_7^3 Y_3 + 2 \omega_1^4 Y_4 + \omega_1^5 Y_5), \nonumber \\
dY_6  +  \omega_1^1 Y_6 & = & \rho^2 (2 \omega_7^3 Y_3 +  \omega_1^4 Y_4 - \omega_1^5 Y_5), \\
dY_3 & = & \omega_7^3 Z_3 + \rho (\omega_1^5 Y_4 + \omega_1^4 Y_5), \nonumber \\
dY_4 & = & - \omega_1^4 Z_4 + \rho (- \omega_1^5 Y_3 + \omega_7^3 Y_5), \nonumber \\
dY_5 & = & \omega_1^5 Z_5 + \rho (- \omega_1^4 Y_3 - \omega_7^3 Y_4), \nonumber 
\end{eqnarray}
where
\begin{eqnarray}
\label{eq:5.66}
Z_3 & = & - Y_6 + \rho^2 (- Y_1 -2 Y_7),\nonumber \\
Z_4 & = & Y_2 + \rho^2 (2 Y_1 + Y_7), \\
Z_5 & = & - Y_2 + Y_6 + \rho^2 (- Y_1 + Y_7). \nonumber 
\end{eqnarray}
From this, we notice that
\begin{equation}
\label{eq:5.67}
Z_3 + Z_4 + Z_5 = 0,
\end{equation}
so that the points $Z_3$, $Z_4$ and $Z_5$ lie on a line.  

From (\ref{eq:5.20}) and (\ref{eq:5.22}) and the lemma, we see that
\begin{equation}
\label{eq:5.68}
d \rho + \rho \omega_1^1 = 0, \quad d \omega_1^1 = 0.
\end{equation}

We now make a change of frame of the form
\begin{eqnarray}
\label{eq:5.69}
Y_1^* & = & \rho Y_1, \quad Y_2^* =  (1/\rho) Y_2, \nonumber \\
Y_7^* & = & \rho Y_7, \quad Y_6^* =  (1/\rho)Y_6, \\
Y_i^* & = & Y_i, \quad 3 \leq i \leq 5. \nonumber 
\end{eqnarray}
Then set
\begin{eqnarray}
\label{eq:5.70}
Z_i^* & = & (1/\rho) Z_i, \quad 3 \leq i \leq 5, \\
\omega_1^{4*} & = & \rho \omega_1^4, \quad  \omega_1^{5*} =  \rho \omega_1^5,
\quad \omega_7^{3*} =  \rho \omega_7^3.\nonumber 
\end{eqnarray}
The effect of this change is to make $\rho^* = 1$ and $\omega_1^{1*} = 0$, for we can compute the following:
\begin{eqnarray}
\label{eq:5.71}
dY_1^* & = & \omega_1^{4*} Y_4 + \omega_1^{5*} Y_5, \nonumber \\
dY_7^* & = & \omega_7^{3*} Y_3 - \omega_1^{5*} Y_5, \nonumber \\
dY_2^* & = & \omega_7^{3*} Y_3 + 2 \omega_1^{4*} Y_4 + \omega_1^{5*} Y_5, \nonumber \\
dY_6^* & = & 2 \omega_7^{3*} Y_3 +  \omega_1^{4*} Y_4 - \omega_1^{5*} Y_5, \\
dY_3 & = & \omega_7^{3*} Z_3^* + \omega_1^{5*} Y_4 + \omega_1^{4*} Y_5, \nonumber \\
dY_4 & = & - \omega_1^{4*} Z_4^* - \omega_1^{5*} Y_3 + \omega_7^{3*} Y_5, \nonumber \\
dY_5 & = & \omega_1^{5*} Z_5^* - \omega_1^{4*} Y_3 - \omega_7^{3*} Y_4, \nonumber 
\end{eqnarray}
with
\begin{eqnarray}
\label{eq:5.72}
dZ_3^* & = & 2 (-2 \omega_7^{3*} Y_3 - \omega_1^{4*} Y_4 + \omega_1^{5*} Y_5),\nonumber \\
dZ_4^* & = & 2 (\omega_7^{3*} Y_3 + 2 \omega_1^{4*} Y_4 + \omega_1^{5*} Y_5), \\
dZ_5^* & = & 2 (\omega_7^{3*} Y_3 - \omega_1^{4*} Y_4 -2 \omega_1^{5*} Y_5),\nonumber 
\end{eqnarray}
and 
\begin{eqnarray}
d\omega_1^{4*} & = & - \omega_1^{5*} \wedge \omega_7^{3*}, \quad {\mbox {\rm i.e., }}\quad
d\theta_1^* = - \theta_2^* \wedge \theta_3^*, \nonumber \\
d\omega_1^{5*} & = & - \omega_7^{3*} \wedge \omega_1^{4*}, \quad {\mbox {\rm i.e., }}\quad
d\theta_2^* = - \theta_3^* \wedge \theta_1^*,  \nonumber \\
d\omega_7^{3*} & = & - \omega_1^{4*} \wedge \omega_1^{5*}, \quad {\mbox {\rm i.e., }}\quad
d\theta_3^* = - \theta_1^* \wedge \theta_2^*. \nonumber 
\end{eqnarray}
Comparing the last equation with (\ref{eq:5.13}), we see that 
$\omega_1^{1*} = 0$ and $\rho^* = 1$.

This is the final frame which we will need in this case $\rho \neq 0$. 
So, we again drop the asterisks.

We are now ready to prove Pinkall's classification result for the case $\rho \neq 0$
\cite[p. 117]{P3}. 
As with the  cyclides, there is only 
one compact model, up to Lie equivalence.  This is  Cartan's isoparametric hypersurface $M^3$ in $S^4$ (see \cite{Car2}--\cite{Car5}).  
It is a tube of constant radius over each of its two focal submanifolds, which
are standard Veronese surfaces in $S^4$.  (See \cite[pp. 296--299]{CR7},
\cite[pp. 151--155 ]{CR8},  \cite{Th6} for more detail.)  We will describe the Veronese surface after stating the theorem.\\

\noindent
{\bf Theorem 5.2:} (Pinkall \cite{P3})  {\it (a) Every connected Dupin submanifold with $\rho \neq 0$
is contained in a unique compact connected Dupin submanifold with $\rho \neq 0$.}\\
{\it (b) Any two Dupin submanifolds with $\rho \neq 0$ are locally Lie equivalent, each
being Lie equivalent to an open subset of Cartan's isoparametric hypersurface in $S^4$.}\\

Our method of proof differs from that of Pinkall in that we will prove directly that each of the focal submanifolds can naturally be considered to be an open subset of a Veronese surface in a hyperplane $P^5 \subset P^6$.  The Dupin submanifold can then be constructed from these focal submanifolds.

We now recall the definition of a Veronese surface.  First consider the map from the unit sphere $y_1^2 + y_2^2 + y_3^2 = 1$ in ${\bf R}^3$ into ${\bf R}^5$ given by
\begin{equation}
\label{eq:5.73}
(x_1,\ldots,x_5) = (2 y_2 y_3, 2 y_3 y_1, 2 y_1 y_2, y_1^2, y_2^2).
\end{equation}
This map takes the same value on antipodal points of the 2-sphere, so it induces a map 
$\phi:P^2 \rightarrow {\bf R}^5$.  One can show by an elementary direct calculation that $\phi$ is an embedding of $P^2$ and that $\phi$ is substantial in ${\bf R}^5$, i.e., does not lie in any hyperplane.  Any embedding of $P^2$ into $P^5$ which is projectively equivalent to $\phi$ is called a {\em Veronese surface}.  (See Lane \cite[pp. 424--430]{Lane} for more detail.)

Let $k_1 = Y_1$, $k_2 = Y_7$, $k_3 = Y_1 + Y_7$ be the focal point maps of the Dupin submanifold $\lambda:B \rightarrow \Lambda$ with $\rho \neq 0$.  Each $k_i$ is constant along the leaves of its principal foliation $T_i$, so each $k_i$ factors through an immersion $\phi_i$ of the 2-dimensional space of leaves $B/T_i$ into $P^6$.  We will show that each of these $\phi_i$ is an open subset of a Veronese surface in some $P^5 \subset P^6$.

We wish to integrate the differential system (\ref{eq:5.71}), which is completely integrable.  For this purpose we drop the asterisks and write the system as follows:
\begin{eqnarray}
\label{eq:5.74}
dY_1 & = & \theta_1 Y_4 + \theta_2 Y_5, \nonumber \\
dY_7 & = & \theta_3 Y_3 - \theta_2 Y_5, \nonumber \\
dY_2 & = & \theta_3 Y_3 + 2 \theta_1 Y_4 + \theta_2 Y_5, \nonumber \\
dY_6 & = & 2 \theta_3 Y_3 +  \theta_1 Y_4 - \theta_2 Y_5, \\
dY_3 & = & \theta_3 Z_3 + \theta_2 Y_4 + \theta_1 Y_5, \nonumber \\
dY_4 & = & - \theta_1 Z_4 - \theta_2 Y_3 + \theta_3 Y_5, \nonumber \\
dY_5 & = & \theta_2 Z_5 - \theta_1 Y_3 - \theta_3 Y_4, \nonumber 
\end{eqnarray}
with
\begin{eqnarray}
\label{eq:5.75}
dZ_3 & = & 2 (-2 \theta_3 Y_3 - \theta_1 Y_4 + \theta_2 Y_5),\nonumber \\
dZ_4 & = & 2 (\theta_3 Y_3 + 2 \theta_1 Y_4 + \theta_2 Y_5), \\
dZ_5 & = & 2 (\theta_3 Y_3 - \theta_1 Y_4 - 2 \theta_2 Y_5),\nonumber 
\end{eqnarray}
where
\begin{eqnarray}
\label{eq:5.76}
d\theta_1 & = & - \theta_2 \wedge \theta_3, \nonumber \\
d\theta_2 & = & - \theta_3 \wedge \theta_1,  \\
d\theta_3 & = & - \theta_1 \wedge \theta_2. \nonumber 
\end{eqnarray}
and
\begin{eqnarray}
\label{eq:5.77}
Z_3 & = & - Y_1 - Y_6 - 2Y_7,\nonumber \\
Z_4 & = & 2Y_1 + Y_2 + Y_7, \\
Z_5 & = & -Y_1 - Y_2 + Y_6 + Y_7, \nonumber 
\end{eqnarray}
so that
\begin{equation}
\label{eq:5.78}
Z_3 + Z_4 + Z_5 = 0.
\end{equation}
Put
\begin{equation}
\label{eq:5.79}
W_1 = - Y_1 + Y_6 -2Y_7, \quad W_2 = -2 Y_1 + Y_2 - Y_7.
\end{equation}
We find from (\ref{eq:5.74}) that
\begin{equation}
\label{eq:5.80}
dW_1=dW_2=0,
\end{equation}
so that the points $W_1$ and $W_2$ are fixed.  Their inner products are
\begin{equation}
\label{eq:5.81}
\langle W_1, W_1 \rangle = \langle W_2, W_2 \rangle =  - 4, \quad \langle W_1, W_2 \rangle = -2,
\end{equation}
and the line $[W_1,W_2]$ consists entirely of timelike points. Its orthogonal complement in 
${\bf R}_2^7$ is spanned by $Y_3,Y_4,Y_5,Z_4,Z_5$.  It consists entirely of spacelike points and
has no point in common with $Q^5$. We will denote it as ${\bf R}^5 $.

It suffices to solve the system (\ref{eq:5.74}) in ${\bf R}^5$ for $Y_3,Y_4,Y_5,Z_4,Z_5$.  For we have 
\begin{equation}
\label{eq:5.82}
d(Z_4-Z_5-6 Y_1)= 0, \quad d(Z_4 + 2Z_5 - 6Y_7) = 0,
\end{equation}
so that there exist constant vectors $C_1$ and $C_2$ such that
\begin{equation}
\label{eq:5.83}
Z_4-Z_5-6 Y_1 = C_1, \quad Z_4 + 2Z_5 - 6Y_7 = C_2.
\end{equation}
Thus, $Y_1$ and $Y_7$ are determined by these equations, and then $Y_2$ and $Y_6$ are determined by  (\ref{eq:5.79}).  Note that $C_1$ and $C_2$ are timelike points, and the line 
$[C_1,C_2]$ consists entirely of timelike points.

Equations (\ref{eq:5.76}) are the  structure equations of $SO(3)$.
It is thus natural to take the latter as the parameter space, whose points are the $3 \times 3$  matrices 
\begin{displaymath}
A = [a_{ik}], \quad 1 \leq i,j,k \leq 3,
\end{displaymath}
satisfying
\begin{equation}
\label{eq:5.84}
A^t A = A A^t = I, \quad \det A = 1.
\end{equation}
The first equations above, when expanded, are
\begin{equation}
\label{eq:5.85}
\sum a_{ij}a_{ik} = \sum a_{ji}a_{ki} = \delta_{jk}.
\end{equation}
The Maurer--Cartan forms of $SO(3)$ are
\begin{equation}
\label{eq:5.86}
\alpha_{ik} = \sum a_{kj} da_{ij} = - \alpha_{ki}.
\end{equation}
They satisfy the Maurer--Cartan equations,
\begin{equation}
\label{eq:5.87}
d \alpha_{ik} = \sum \alpha_{ij} \wedge \alpha_{jk}.
\end{equation}
If we set
\begin{equation}
\label{eq:5.88}
\theta_1 = \alpha_{23}, \quad \theta_2 = \alpha_{31}, \quad \theta_3 = \alpha_{12},
\end{equation}
these equations reduce to (\ref{eq:5.76}).  With the $\theta_i$ given by (\ref{eq:5.88}),
we shall write down an explicit solution of (\ref{eq:5.74}).

Let $E_A, 1 \leq A \leq 5$, be a fixed linear frame in ${\bf R}^5$.  Let
\begin{equation}
\label{eq:5.89}
F_i = 2a_{i2}a_{i3}E_1 + 2a_{i3}a_{i1}E_2 + 2a_{i1}a_{i2}E_3 + a_{i1}^2 E_4 + a_{i2}^2 E_5, 
\end{equation}
$1 \leq i,j,k \leq 3$. Since
\begin{displaymath}
a_{i1}^2 + a_{i2}^2 + a_{i3}^2 = 1,
\end{displaymath}
we see from (\ref{eq:5.73}), with $y_j = a_{ij}$, that $F_i$ is a Veronese surface for $1 \leq i \leq 3$.  Using equation (\ref{eq:5.85}), we compute that
\begin{equation}
\label{eq:5.90}
F_1 + F_2 + F_3 = E_4 + E_5 = {\mbox {\rm constant}}.
\end{equation}
Since the coefficients in $F_i$ are quadratic, the partial derivatives,
\begin{displaymath}
\frac{\partial^2 F_i}{\partial a_{ij} \partial a_{ik}},
\end{displaymath}
are independent of $i$.  Moreover, the quantities $G_{ik}$ defined below satisfy
\begin{displaymath}
G_{ik} = \sum_j a_{ij} \frac{\partial F_k}{\partial a_{kj}} = G_{ki}.
\end{displaymath}
We use these facts in the following computation:
\begin{eqnarray}
\label{eq:5.92}
dG_{ik} & = & \sum \frac{\partial F_k}{\partial a_{kj}} da_{ij} + 
\sum a_{ij} \frac{\partial^2 F_k}{\partial a_{kj} \partial a_{kl}}da_{kl}\\
& = & \sum \frac{\partial F_k}{\partial a_{kj}} da_{ij} +
\sum a_{ij} \frac{\partial^2 F_i}{\partial a_{ij} \partial a_{il}}da_{kl}\nonumber \\
& = & \sum \frac{\partial F_k}{\partial a_{kj}} da_{ij} +
\sum \frac{\partial F_i}{\partial a_{ij}} da_{kj},\nonumber
\end{eqnarray}
where the last step follows from the linear homogeneity of $\partial F_i/\partial a_{il}$.  In terms of $\alpha_{ij}$, we have
\begin{equation}
\label{eq:5.93}
dG_{ik} = \sum \frac{\partial F_k}{\partial a_{kj}}a_{lj}\alpha_{il}
+ \sum \frac{\partial F_i}{\partial a_{ij}}a_{lj}\alpha_{kl},
\end{equation}
which gives, when expanded,
\begin{equation}
\label{eq:5.94}
dG_{23} = 2(F_3 - F_2) \theta_1 + G_{12} \theta_2 - G_{13} \theta_3,
\end{equation}
and its cyclic permutations.

On the other hand, by the same manipulation, we have
\begin{equation}
\label{eq:5.95}
dF_i = \sum \frac{\partial F_i}{\partial a_{ij}} da_{ij} = 
\sum \frac{\partial F_i}{\partial a_{ij}} a_{kj} \alpha_{ik},
\end{equation}
giving
\begin{displaymath}
dF_1 = - G_{31}\theta_2 + G_{12} \theta_3,
\end{displaymath}
and its cyclic permutations.

One can now immediately verify that a solution of (\ref{eq:5.74}) is given by
\begin{eqnarray}
\label{eq:5.96}
Y_3 = G_{12}, \quad Y_4 & = & - G_{23}, \quad Y_5 = G_{31}, \\
Z_3 = 2(F_2-F_1), \quad Z_4 & = & 2(F_3-F_2), \quad Z_5 = 2(F_1-F_3).\nonumber
\end{eqnarray}
This is also the most general solution of (\ref{eq:5.74}), for the solution is determined up to a linear transformation, and our choice of frame $E_A$ is arbitrary.

By (\ref{eq:5.96}), the functions $Z_1$, $Z_2$, $Z_3$ are expressible in terms of $F_1$, $F_2$, $F_3$, and then by (\ref{eq:5.83}), so also are $Y_1$, $Y_7$, $Y_1+Y_7$.  Specifically, by 
(\ref{eq:5.83}), (\ref{eq:5.90}) and (\ref{eq:5.96}), we have
\begin{eqnarray}
6Y_1 & = & Z_4 - Z_5 - C_1 = 2(-F_1 -F_2 + 2F_3) - C_1 \nonumber \\
 & = & 6F_3 - 2(E_4+E_5) - C_1, \nonumber
\end{eqnarray}
so that the focal map $Y_1$, up to an additive constant, is the Veronese surface $F_3$.  Similarly,
the focal maps $Y_7$ and $Y_1+Y_7$ are the Veronese surfaces $F_1$ and $-F_2$, respectively, up to additive constants.

We see from (\ref{eq:5.79}) that
\begin{displaymath}
\langle Y_1, W_1 \rangle = 0, \quad \langle Y_7, W_2 \rangle = 0, \quad \langle Y_1+Y_7, W_1-W_2 \rangle = 0.
\end{displaymath}
Thus $Y_1$ is contained in the M\"{o}bius space $\Sigma^4 = Q^5 \cap W_1^{\perp}$. Let
$e_7 = W_1/2$ and 
$e_1 = (2W_2 - W_1)/\sqrt{12}$. 
Then $e_1$ is the unique unit vector on the timelike line $[W_1,W_2]$
which is orthogonal to $W_1$.  In a manner similar to that of Section 3, we can write
\begin{displaymath}
Y_1 = e_1 + f,
\end{displaymath}
where $f$ maps $B$ into the unit sphere $S^4$ in the Euclidean space 
${\bf R}^5 = [e_1,e_7]^{\perp}$ in ${\bf R}^7_2$.
We call $f$ the {\em spherical projection} of the Legendre map $\lambda$ determined by the
ordered pair $\{e_7, e_1\}$ (see \cite{CC1} for more detail).  
We know that $f$ is constant along the leaves of the principal foliation $T_1$ 
corresponding to $Y_1$, and $f$
induces a map $\tilde{f}: B/T_1 \rightarrow S^4$.
By what we have shown above, $\tilde{f}$ must be an open subset of a spherical Veronese surface.

Note that the unit timelike vector $W_2/2$ satisfies
\begin{displaymath}
W_2/2 = (\sqrt{3}/2)e_1 + (1/2) e_7 = \sin(\pi/3) e_1 + \cos(\pi/3) e_7.
\end{displaymath}
If we consider the points in the M\"{o}bius space $\Sigma$ to represent point spheres in $S^4$, then as we show in \cite{CC1}, points in $Q^5 \cap W_2^{\perp}$ represent oriented
spheres in $S^4$ with oriented radius $-\pi/3$.  In a way similar to that above, the second focal submanifold $Y_7 \subset Q^5 \cap W_2^{\perp}$ induces a spherical Veronese surface. 
When considered from the point of view of the M\"{o}bius space $\Sigma$, the points $Y_7$ represent  oriented spheres of radius
$-\pi/3$ centered at points of this Veronese surface.  These spheres must be in
oriented contact with the point spheres of the first Veronese surface determined by $Y_1$ in 
$Q^5 \cap {\bf P}^5$.
Thus, the points in the second Veronese surface must lie at a distance
$\pi/3$ along normal geodesics in $S^4$ to the 
first Veronese surface $\tilde{f}$.  In fact (see, for example \cite[pp. 296--299]{CR7}), the set of all points in $S^4$ at a distance $\pi/3$ from a spherical Veronese surface is another spherical Veronese surface.

Thus, with this choice of coordinates, the Dupin submanifold in question is simply an open subset of the Legendre  submanifold induced in the standard way by considering $B$
to be the unit normal bundle to the spherical Veronese embedding
$\tilde{f}$ induced by $Y_1$.  For values
of $t = k\pi/3, k \in {\bf Z}$, the parallel hypersurface
at oriented distance $t$ to $\tilde{f}$ in $S^4$ is a
Veronese surface.  For other values of $t$, the parallel hypersurface is an isoparametric hypersurface in $S^4$ with three distinct principal curvatures (Cartan's isoparametric hypersurface).
All of these parallel hypersurfaces are 
Lie equivalent to each other and to the Legendre submanifolds induced by the Veronese surfaces.\\

\noindent
{\bf Remark 5.2.} (added by T. Cecil in 2020) There is a somewhat shorter way to complete the proof of Theorem 5.2 by using the following characterization of isoparametric hypersurfaces in Lie sphere geometry (see \cite{Cec4}, \cite[p. 77]{Cec1} or  \cite[p. 221]{CR8}).  
Let  $\lambda: M^{n-1} \rightarrow \Lambda^{2n-1}$ be a Legendre submanifold with $g$ distinct curvature spheres $K_1,\ldots, K_g$ at each point.  Then $\lambda$ is Lie equivalent 
to the Legendre lift of an isoparametric hypersurface with $g$ principal curvatures in $S^n$ if and only if there exist $g$ points $P_1,\ldots,P_g$, on a timelike line in $P^{n+2}$ such that $\langle K_i,P_i \rangle = 0,$  $1 \leq i \leq g$.  
In the proof of Theorem 5.2 above, we see that the three curvature sphere maps
$Y_1$, $Y_7$ and  $Y_1 + Y_7$ are orthogonal to the three points  
$W_1$, $W_2$, and $W_1 - W_2$ on the 
timelike line in $[W_1,W_2]$ in $P^6$ by equations (\ref{eq:5.79})--(\ref{eq:5.81}).  Thus,
$\lambda$ is Lie equivalent to the Legendre lift of an isoparametric hypersurface with three principal curvatures in $S^4$.  The proof of Theorem 5.2 is then completed by invoking Cartan's \cite{Car3} classification of such isoparametric hypersurfaces.\\

\noindent
{\em The case $\rho = 0$.}\\

\noindent
We now consider the case where $\rho$ is identically zero on $B$.  It turns out that no new examples occur here, in that these Dupin submanifolds can all be constructed from Dupin cyclides by certain standard constructions.  To make this precise, we recall's Pinkall's \cite[p. 437]{P4} 
notion of reducibility.  Our Dupin submanifold can be considered, as in Section 3, to have been induced from a Dupin hypersurface $M^3 \subset E^4$.  The Dupin submanifold is {\em reducible} if $M^3$ is obtained from a Dupin surface $S \subset E^3 \subset E^4$ by one of the four following standard constructions.
\begin{eqnarray}
\label{eq:5.97}
&i.& M\  {\rm is\ a\ cylinder}\ S \times {\bf R} \ {\rm in}\ E^4.           \nonumber \\
&ii.& M\ {\rm is\ the\ hypersurface\ of\ revolution\ obtained\ by\ revolving}\ \nonumber \\
&   & S\ {\rm about\  a\  plane}\ \pi \ {\rm disjoint\ from}\ S \ {\rm in}\ E^3. \\
&iii.& {\rm Project}\ S\ {\rm stereographically\ onto\ a\ surface}\ N \subset S^3 \subset E^4. \nonumber\\
&   & M \ {\rm is\ the\ cone}\  {\bf R} \cdot N\ {\rm over}\ N.\nonumber \\
&iv.& M\ {\rm is\ a\ tube\ of\ constant\ radius\ around}\ S \ {\rm in}\ E^4.\nonumber
\end{eqnarray}

\noindent
{\bf Remark 5.3.} (added by T. Cecil in 2020) 
More generally, the Dupin submanifold $\lambda:B \rightarrow \Lambda^7$
is said to be reducible
if it is locally Lie equivalent to the Legendre lift of a proper Dupin hypersurface in $E^4$ obtained by one of Pinkall's standard constructions
given in equation (\ref{eq:5.97}).  See \cite[p. 437]{P4}  or \cite[pp.  125--148]{Cec1} for more detail on Pinkall's standard constructions in the context of Lie sphere geometry.\\

\noindent
Pinkall proved \cite[p. 438]{P4} that the Dupin submanifold $\lambda:B \rightarrow \Lambda^7$ is reducible if and only if some focal point map is contained in a 4-dimensional subspace 
$P^4 \subset P^6$.

If $\rho$ is identically zero on $B$, then by  (\ref{eq:5.20}), all of the covariant derivatives of $\rho$ are also equal to zero.  From  (\ref{eq:5.21}) and  (\ref{eq:5.54}), we see that the functions in equations  (\ref{eq:5.10})--(\ref{eq:5.12}) satisfy
\begin{displaymath}
q=t=b=0, \quad a=p=0,\quad r=s, \quad c=-u.
\end{displaymath}
Then from (\ref{eq:5.55}), we have $\omega_2^7 = 0$.  From these and the other relations among the Maurer-Cartan forms which we have derived, we see that the differentials of the frame vectors can be written
\begin{eqnarray}
\label{eq:5.98}
dY_1  -  \omega_1^1 Y_1 & = & \omega_1^4 Y_4 + \omega_1^5 Y_5, \nonumber \\
dY_7  -  \omega_1^1 Y_7 & = & \omega_7^3 Y_3 - \omega_1^5 Y_5, \nonumber \\
dY_2  +  \omega_1^1 Y_2&  = & s (- \omega_1^4 Y_4 + \omega_1^5 Y_5), \nonumber \\
dY_6  +  \omega_1^1 Y_6 & = & u (\omega_7^3 Y_3 + \omega_1^5 Y_5), \\
dY_3 & = & \omega_7^3 (-Y_6 + u Y_7), \nonumber \\
dY_4 & = & \omega_1^4 (sY_1 - Y_2), \nonumber \\
dY_5 & = & \omega_1^5 (-sY_1 - Y_2 + Y_6 -u Y_7). \nonumber 
\end{eqnarray}
Note that from (\ref{eq:5.44}) and (\ref{eq:5.45}), we have
\begin{equation}
\label{eq:5.99}
ds + 2s \omega_1^1 = 0, \quad du + 2u \omega_1^1 = 0,
\end{equation}
and from (\ref{eq:5.13}) that
\begin{equation}
\label{eq:CC-5.100}
d \theta_i = \omega_1^1 \wedge \theta_i, \quad i = 1,2,3.
\end{equation}
From  (\ref{eq:5.22}), we have that $d\omega_1^1 = 0$.  
Hence on any local disk neighborhood $U$ in $B$, we have that
\begin{equation}
\label{eq:5.101}
\omega_1^1 = d\sigma,
\end{equation}
for some smooth function
$\sigma$ on $U$. We next consider a change of frame of the form,
\begin{eqnarray}
\label{eq:5.102}
Y_1^* & = & e^{-\sigma} Y_1, \quad Y_7^* =  e^{-\sigma} Y_7, \nonumber \\
Y_2^* & = & e^{\sigma} Y_2, \quad Y_6^* =  e^{\sigma} Y_6,\\
Y_i^* & = & Y_i, \quad i = 3,4,5. \nonumber 
\end{eqnarray}
The effect of this change is to make $\omega_1^{1*} = 0$ while keeping $\rho^* = 0$.
If we set
\begin{displaymath}
\omega_1^{4*} = e^{- \sigma} \omega_1^4, \quad \omega_1^{5*} = e^{- \sigma} \omega_1^5, \quad
\omega_7^{3*} = e^{- \sigma} \omega_7^3,  
\end{displaymath}
then we can compute that from (\ref{eq:5.98}) that
\begin{eqnarray}
\label{eq:5.103}
dY_1^* & = & \omega_1^{4*} Y_4 + \omega_1^{5*} Y_5, \nonumber \\
dY_7^* & = & \omega_7^{3*} Y_3 - \omega_1^{5*} Y_5, \nonumber \\
dY_2^* & = & s^*(- \omega_1^{4*} Y_4 + \omega_1^{5*} Y_5) \nonumber \\
dY_6^* & = & u^*(\omega_7^{3*} Y_3 + \omega_1^{5*} Y_5),\\
dY_3 & = & \omega_7^{3*} Z_3^*,\ {\rm where}\   Z_3^* =  - Y_6^* - u^* Y_7^*,\nonumber \\
dY_4 & = &  \omega_1^{4*} Z_4^*,\   {\rm where}\ Z_4^* =  s^*Y_1^* - Y_2^*,\nonumber \\
dY_5 & = & \omega_1^{5*} Z_5^*, \  {\rm where}\ Z_5^* =  -s^* Y_1^* - Y_2^* + Y_6^* - u^* Y_7^*,  \nonumber 
\end{eqnarray}
where
\begin{equation}
\label{eq:5.104}
s^* = s e^{2\sigma}, \quad u^* = u e^{2\sigma}.
\end{equation}
Using (\ref{eq:5.99}) and (\ref{eq:5.104}), we can then compute that
\begin{equation}
\label{eq:5.105}
ds^* = 0, \quad du^* = 0,
\end{equation}
i.e., $s^*$ and $u^*$ are constant functions on the local neighborhood $U$.

The frame in equation (\ref{eq:5.102}) is our final frame, and we drop the asterisks
in further references to equations (\ref{eq:5.102})--(\ref{eq:5.105}).  Since the functions
$s$ and $u$ are now constant, we can compute from equation (\ref{eq:5.103}) that
\begin{eqnarray}
\label{eq:5.106}
dZ_3 & = & -2u \omega_7^3 Y_3,\nonumber \\
dZ_4 & = & 2s \omega_1^4 Y_4, \\
dZ_5 & = & 2(u-s) \omega_1^5 Y_5. \nonumber 
\end{eqnarray}
From this we see that the following $4$-dimensional subspaces,
\begin{eqnarray}
\label{eq:5.107}
&{\mbox {\rm Span}}&\{Y_1,Y_4,Y_5,Z_4,Z_5\},\nonumber \\
&{\mbox {\rm Span}}&\{Y_7,Y_3,Y_5,Z_3,Z_5\}, \\
&{\mbox {\rm Span}}&\{Y_1+Y_7,Y_3,Y_4,Z_3,Z_4\},\nonumber 
\end{eqnarray}
are invariant under exterior differentiation and are thus constant.  Thus, each of the three focal point maps $Y_1$, $Y_7$ and $Y_1+Y_7$ is contained in a $4$-dimensional subspace of $P^6$, and
our Dupin submanifold is reducible in three different ways.  Each of the three focal point maps is thus an immersion of the space of leaves of its principal foliation onto an open subset of a cyclide of Dupin in a space $\Sigma^3 = P^4 \cap Q^5$.  See \cite[p. 188]{Cec1} for more detail.\\

We state this result to due Pinkall \cite{P3} as follows:\\

\noindent
{\bf Theorem 5.3:} Every Dupin submanifold with $\rho = 0$ is reducible.  Thus, it is obtained from a cyclide in ${\bf R}^3$ by one of the four standard constructions (\ref{eq:5.97}).\\

Pinkall \cite[p. 111]{P3} then proceeds to classify Dupin submanifolds with $\rho = 0$ up to Lie equivalence.  We will not prove his result here.  The reader can follow his proof using the fact that his constants $\alpha$ and $\beta$ are our constants $s$ and $-u$, respectively.\\

\noindent
{\bf Remark 5.4.} (added by T. Cecil in 2020) 
We now briefly describe some local results for higher dimensional Dupin hypersurfaces that have been obtained using an approach similar to that used in proving Theorems 5.2 and 5.3 above.
We also describe some remaining open problems in the field. (See \cite[pp.  301--308]{CR8}, \cite{Cec9}, for more complete surveys.)

These classifications primarily involve finding proper Dupin submanifolds that are irreducible in the sense of Pinkall \cite[p. 438]{P4} as in Theorem 5.2. In general, a proper Dupin submanifold
$\lambda: M^{n-1} \rightarrow \Lambda^{2n-1}$ is said to be {\em reducible}
if it is locally Lie equivalent to the Legendre lift of a proper Dupin hypersurface in ${\bf R}^n$ obtained by one of Pinkall's standard constructions (suitably generalized to ${\bf R}^n$)
given in equation (\ref{eq:5.97}).  

Pinkall \cite{P4} proved the following characterization of reducibility in terms of Lie sphere geometry.
A connected proper Dupin submanifold $\lambda: M^{n-1} \rightarrow \Lambda^{2n-1}$ 
is reducible if and only if there exists a
curvature sphere map $[K]$ of $\lambda$ that lies in a linear subspace
of $P^{n+2}$ of codimension at least two. (See also \cite[pp. 127--148]{Cec1} and \cite[pp. 233--256]{CR8} for more detail.)

After Pinkall's results, Niebergall \cite{N1} next proved that every connected proper Dupin hypersurface in ${\bf R}^5$ with three principal curvatures is reducible. 
Following that, Cecil and Jensen \cite{CJ2} found the following local classification of irreducible Dupin hypersurfaces with three principal curvatures of arbitrary dimension.
Suppose $M^{n-1}$ is a connected irreducible
proper Dupin hypersurface in $S^n$ with
three distinct principal curvatures of multiplicities
$m_1,m_2,m_3$.  Then $m_1 = m_2 = m_3 = m$, and $M^{n-1}$ is Lie 
equivalent to the Legendre lift of an isoparametric hypersurface in $S^n$.   
 
It then follows from Cartan's \cite{Car3} classification of
isoparametric hypersurfaces with $g=3$  that $m=1,2,4$ or 8, and the isoparametric hypersurface must be a tube of constant radius
over a standard embedding of a projective
plane ${\bf F}P^2$ into $S^{3m+1}$ (see, for example,  \cite[pp. 296--299]{CR7}), 
 \cite[pp. 151--155]{CR8}), \cite{Th6}),
where ${\bf F}$ is the division algebra
${\bf R}$, ${\bf C}$, ${\bf H}$ (quaternions),
${\bf O}$ (Cayley numbers), for $m=1,2,4,8,$ respectively.  Thus, up to
congruence, there is only one such family for each value of $m$.

Note that in the original paper \cite{CJ2}, the theorem of Cecil and Jensen was proven under the assumption that  $M^{n-1}$ is locally irreducible, i.e., that $M^{n-1}$ does not contain any reducible open subset. However, Cecil, Chi and Jensen \cite{CCJ3} later proved that proper Dupin hypersurfaces are analytic, and as a consequence local irreducibility is equivalent to
irreducibility (see, \cite{CCJ2} or \cite[pp. 145--146]{Cec1}).

An open problem is the classification up to Lie equivalence of reducible Dupin hypersurfaces  with three principal curvatures in ${\bf R}^n, n \geq 5$.
As noted above, Pinkall \cite{P3} found such a classification in the 
case of $M^3 \subset {\bf R}^4$. It may be
possible to generalize Pinkall's results to higher
dimensions using the approach of \cite{CJ2}. 

We next discuss the problem of finding local classifications of irreducible proper Dupin submanifolds $\lambda$ with four curvature spheres.  As mentioned in Remark 1.2, Pinkall and Thorbergsson \cite{PT1} (for $g=4$), and Miyaoka and Ozawa \cite{MO} (for $g=4$ and $g=6$) gave constructions of compact proper Dupin hypersurfaces in spheres that do not have constant Lie curvatures (cross-ratios of the principal curvatures), and so they are not Lie equivalent to the Legendre lift of an isoparametric hypersurface.  These examples are irreducible (see \cite[pp. 255--256]{CR8}).  So irreducibility itself without additional assumptions does not lead to the conclusion that $\lambda$ is Lie equivalent to an isoparametric hypersurface, as in the case $g=3$.

This leads to the problem of finding necessary and sufficient conditions for an irreducible proper Dupin hypersurface with $g=4$ curvature spheres to be Lie equivalent to the Legendre lift of an isoparametric hypersurface.  Considerable progress has been made on this problem, as we now describe.

Let $\lambda:M^{n-1} \rightarrow \Lambda^{2n-1}$
be a connected proper Dupin submanifold with four curvature spheres at each point. Then one can find a Lie frame in which the four curvature spheres are $Y_1$, $Y_{n+3}$, $Y_1+Y_{n+3}$ and $Y_1+ \Psi Y_{n+3}$, 
where $\Psi$ is the Lie curvature of $\lambda$,
having respective multiplicities $m_1$, $m_2$, $m_3$ and $m_4$.
Corresponding to the one function $\rho$ in the case Theorems 5.2 and 5.3, there are four sets of functions,
$F^{\alpha}_{pa},F^{\mu}_{pa},F^{\mu}_{\alpha a},F^{\mu}_{\alpha p}$, 
where
$1  \leq  a \leq m_1$, 
$m_1+1  \leq   p \leq m_1+m_2$,
$m_1+m_2+1  \leq  \alpha \leq m_1+m_2+m_3$,
$m_1+m_2+m_3+1  \leq  \mu \leq n-1$.
 
M\"{u}nzner \cite{Mu}--\cite{Mu2} proved that the multiplicities of the principal curvatures of an isoparametric hypersurface with four principal curvatures in $S^n$ must satisfy $m_1 = m_2, m_3 = m_4$, when the principal curvatures are appropriately ordered (see also \cite{Ab}, \cite{Stolz}, \cite{CCJ1}).  Thus, in the papers of Cecil and Jensen \cite{CJ3}, and Cecil, Chi and Jensen
\cite{CCJ2}, those restrictions on the multiplicities are assumed.  Furthermore, the assumption that 
that the Lie curvature $\Psi = -1$ is made in \cite{CCJ2}, since that is also
is true for an isoparametric hypersurface with four principal curvatures (with the appropriate ordering  of the principal curvatures).  

In \cite[pp. 3--4]{CJ3}, Cecil and Jensen conjectured that an irreducible 
connected proper Dupin hypersurface in 
$S^n$ with four principal curvatures having multiplicities satisfying $m_1 = m_2, m_3 = m_4$ and constant Lie curvature
$\Psi$ must be Lie equivalent to the Legendre lift of an open subset of an isoparametric hypersurface in $S^n$.

In that paper \cite{CJ3}, Cecil and Jensen proved that the conjecture is true if all the multiplicities are equal to one
(see also Niebergall \cite{N2}, who obtained the same conclusion under additional assumptions).  In the second paper \cite{CCJ2} mentioned above, 
Cecil, Chi and Jensen proved that the conjecture is true if $m_1 = m_2 \geq 1$, 
$m_3 = m_4 = 1$, and the Lie curvature is assumed to
satisfy $\Psi = -1$.  
They conjectured that their theorem is still true if
$m_3 = m_4$ is also allowed to be greater than one, but the calculations involved in that case are formidable, and it remains an open problem.

Note that the assumption of irreducibility is needed in the papers \cite{CCJ2}  and \cite{CJ3} mentioned in the preceding paragraph, since
Cecil \cite{Cec4} constructed an example of a reducible noncompact proper Dupin submanifold with $g=4$ curvature spheres having all multiplicities equal and
constant Lie curvature  $\Psi = -1$ that is not Lie equivalent to an open subset of an isoparametric hypersurface with four principal curvatures.
(See also \cite[pp. 80--82]{Cec1} or \cite[pp.304--306]{CR8}.)

An important step in in the theorems of
Cecil-Jensen \cite{CJ3} and Cecil-Chi-Jensen \cite{CCJ2} mentioned above
is proving that under the assumptions $m_1 = m_2 \geq 1$, 
$m_3 = m_4 = 1$, and $\Psi = -1$, the Dupin submanifold $\lambda$ is reducible if there exists some fixed index, say $a$,
such that
\begin{equation}
\label{eq:4.6.87}
F^{\alpha}_{pa} = F^{\mu}_{pa} = F^{\mu}_{\alpha a} = 0, \quad{\mbox {\rm for all }}p, \alpha, \mu.
\end{equation}
Thus, if $\lambda$ is irreducible, no such index $a$ can exist, and it can be shown after a lengthy argument that $\lambda$ is Lie equivalent to the Legendre lift of an open subset
of an isoparametric hypersurface, as 
in Remark 5.2.  This is a generalization of the case $\rho \neq 0$ for $g=3$ handled in Theorem 5.2 above.

\noindent Thomas E. Cecil

\noindent Department of Mathematics and Computer Science

\noindent College of the Holy Cross

\noindent Worcester, MA 01610

\noindent email: tcecil@holycross.edu\\

\noindent
Shiing-Shen Chern

\noindent
Department of Mathematics

\noindent
University of California

\noindent
Berkeley, CA 94720

\noindent
and

\noindent
Mathematical Sciences Research Institute

\noindent
1000 Centennial Drive

\noindent
Berkeley, CA 94720

\end{document}